\documentclass{jgcc}

\pdfoutput=1

\usepackage{lastpage}

\jgccdoi{17}{1}{2}{13954} 

\jgccheading{}{\pageref{LastPage}}{}{}{Jul.~18,~2024}{Apr.~23,~2025}{}   


\usepackage{url}
\usepackage{import}

\usepackage[T1]{fontenc}
\usepackage{lmodern}
\usepackage[english]{babel}
\usepackage{amsmath}
\usepackage{amsfonts}
\usepackage{amssymb}
\usepackage{mathrsfs} 

\usepackage{dsfont}
\usepackage{mathtools} 

\definecolor{TealBlue}{rgb}{0, .68, .7}
\definecolor{Turquoise}{rgb}{0,.7,.81}
\definecolor{LimeGreen}{rgb}{.196,.804,.196}
\definecolor{Magenta}{rgb}{1,0,1}

\usepackage{enumitem}

\usepackage{changepage}
\usepackage{adjustbox}

\usepackage[font=small, hypcap=false]{caption}
\usepackage{parskip}

\usepackage{tikz-cd}
\usepackage{tikz}
\usetikzlibrary{shapes, shapes.geometric, arrows, arrows.meta, decorations.markings, decorations.pathreplacing, automata}
\pgfdeclarelayer{background} 
\pgfdeclarelayer{foreground}
\pgfsetlayers{background,main,foreground}
\newcommand{\vardot}{\,\cdot\,}
\newcommand{\anode}{{\tikz[baseline=-3.5pt] \node[state, accepting, minimum size=8pt, inner sep=0pt] at (0,0) {};}}
\newcommand{\0}{{\bf 0}}

\newcommand{\Conv}{\mathrm{ConvHull}}

\newcommand{\KP}{\mathsf{KS}}
\newcommand{\SMP}{\mathsf{SMM}}
\newcommand{\RSMP}{\mathsf{RatM}}

\newcommand{\start}{\mathsf{start}}
\newcommand{\accept}{\mathsf{accept}}
\DeclareMathOperator{\ev}{ev} 

\newcommand{\say}[1]{``#1''}

\newcommand{\N}{\mathbb N}
\newcommand{\Z}{\mathbb Z}
\newcommand{\Q}{\mathbb Q}
\newcommand{\R}{\mathbb R}

\newcommand{\Sch}{\mathcal Sch}

\newcommand{\Kc}{\mathcal K}
\newcommand{\Lc}{\mathcal L}    
\newcommand{\Zc}{\mathcal Z}    

\newcommand{\reg}{\mathrm{reg}}
\newcommand{\abn}{\mathrm{abn}}

\newcommand{\onto}{\twoheadrightarrow}

\newcommand{\longto}{\longrightarrow}

\newcommand{\la}{\left\langle}
\newcommand{\ra}{\right\rangle}

\newcommand{\abs}[1]{\left| #1 \right|}
\newcommand{\norm}[1]{\left\|#1\right\|}

\renewcommand{\ge}{\geqslant}
\renewcommand{\le}{\leqslant}

\DeclareMathOperator{\rk}{rk}

\usepackage{hyperref}

\newtheorem{innercustomthm}{Theorem}
\newenvironment{customthm}[1]
{\renewcommand\theinnercustomthm{#1}\innercustomthm}
{\endinnercustomthm}

\newtheorem{innercustomprop}{Proposition}
\newenvironment{customprop}[1]
{\renewcommand\theinnercustomprop{#1}\innercustomprop}
{\endinnercustomprop}

\newtheorem{innercustomcor}{Corollary}
\newenvironment{customcor}[1]
{\renewcommand\theinnercustomcor{#1}\innercustomcor}
{\endinnercustomcor}

\theoremstyle{plain}
\newtheorem*{thm*}{Theorem}
\newtheorem*{cor*}{Corollary}

\theoremstyle{definition}
\newtheorem*{rem*}{Remark}

\begin{document}
	
	\title{Membership problems in nilpotent groups}
	
	\author[C.~Bodart]{Corentin Bodart}	
	\address{Mathematical Institute, University of Oxford, UK}	
	\email{cobodart123@gmail.com}  
	\thanks{The author was partially supported by the Swiss NSF grant 200020-200400.}	
	
	
	
	
	
	\begin{abstract}
		We study both the Submonoid Membership problem ($\SMP$) and the Rational Subset Membership problem ($\RSMP$) in finitely generated nilpotent groups. We give two reductions with important applications: 
		\medbreak
		\begin{itemize}[leftmargin=6mm]
			\item The $\SMP$ in any nilpotent group can be reduced to the $\RSMP$ in \emph{smaller} groups. As a corollary, we prove the existence of a group with decidable $\SMP$ and undecidable $\RSMP$, answering a question of Lohrey and Steinberg.
			
			\medbreak
			
			\item The Rational Subset Membership problem in $H_3(\Z)$ can be reduced to the Knapsack problem in the \emph{same} group, and is therefore decidable.
		\end{itemize}
		\medbreak
		\noindent We deduce that the filiform $3$-step nilpotent group has decidable Submonoid Membership.
	
	\end{abstract}
	
	\maketitle
	

\quad Decision problems are central motivating questions in combinatorial group theory. The first example is the Word Problem, introduced by Dehn in 1910 \cite{Dehn1910}. For matrix groups, the first hard problem is the Subgroup Membership problem. With it came one of the first undecidability results in $SL_4(\Z)$ \cite{Mihailova1958}, and also some positive results in nilpotent and polycyclic groups \cite{Malcev}. This paper focuses on the next two decision problems in line, namely the (uniform) Submonoid and Rational Subset Membership problems. For a group $G$ (given as a matrix group, or endowed with a finite generating set $S$) we attempt to produce algorithms with the following specifications:

\textbf{Submonoid Membership problem ($\SMP(G)$)}
\begin{itemize}[leftmargin=18mm, rightmargin=2mm]
	\item[Input:\hspace*{3mm}] Elements $g$ and $g_1,g_2,\ldots,g_n\in G$ (defined as matrices or words over $S$).
	\item[Output:] Decide whether $g$ belongs to the submonoid $\{g_1,g_2,\ldots,g_n\}^*$.
\end{itemize}

\textbf{Rational Subset Membership problem ($\RSMP(G)$)}
\begin{itemize}[leftmargin=18mm, rightmargin=2mm]
	\item[Input:\hspace*{3mm}] An element $g\in G$ (given either as a matrix or as a word over $S$) and a rational subset $R\subseteq G$ (defined by a finite state automaton, labeled by elements in $G$) 
	\item[Output:] Decide whether $g\in R$.
\end{itemize}
Both problems are known to be decidable for a few classes of groups, including free groups \cite{Benois1} and virtually abelian groups \cite{Grunschlag}. We suggest to have a look at the surveys \cite{survey_Lohrey, Dong_Survey, Lohrey_survey24} for a more exhaustive picture. 

Finally, we recall the Knapsack Problem which plays a key role in our paper. This problem is a special case of Rational Subset Membership (see Figure \ref{fig:knapsack}), introduced in \cite{Knapsack}.

\textbf{Knapsack problem ($\KP(G)$)}
\begin{itemize}[leftmargin=18mm, rightmargin=2mm]
	\item[Input:\hspace*{3mm}] Elements $g$ and $g_1,g_2,\ldots,g_n\in G$ (defined as matrices or words over $S$).
	\item[Output:] Decide whether $g\in \{g_1\}^*\,\{g_2\}^*\ldots \{g_n\}^*$.
\end{itemize}
This sub-problem is known to be decidable in a much larger class of group. This includes hyperbolic groups \cite{Knapsack} and co-context-free groups \cite{Knapsack_Heis}. 
\begin{center}
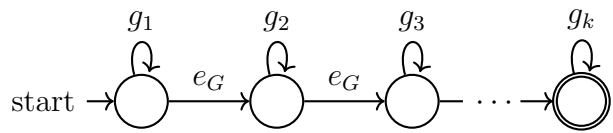

	\begin{tikzpicture}[scale=.45, thick]	
		\node[state, initial left, initial distance=8mm, minimum size=20pt] (v0) at (0,0) {};
		\node[state, minimum size=20pt] (v1) at (4,0) {};
		\node[state, minimum size=20pt] (v2) at (8,0) {};
		\node[state, accepting, minimum size=20pt,] (v3) at (13,0) {};
		
		\path (v0)  edge[->] node [above] {$e_G$} (v1)
		(v0)  edge [loop above] node [above] {$g_1$} (v0)
		(v1)  edge[->] node [above] {$e_G$} (v2)
		(v1)  edge [loop above] node [above] {$g_2$} (v1)
		(v2)  edge [loop above] node [above] {$g_3$} (v2)
		(v2)  edge[-] (9.5,0)
		(11.2,0)  edge[->] (v3)
		(v3)  edge [loop above] node [above] {$g_k$} (v3);
		
		\node at (10.5,-.05) {$\cdots$};
	\end{tikzpicture}
	\captionsetup{margin=10mm, font=footnotesize}
	\captionof{figure}{Finite state automaton relative to the Knapsack problem}
	\label{fig:knapsack}
\end{center}

Finitely generated submonoids are rational, hence the decidability of $\RSMP(G)$ implies the decidability of $\SMP(G)$. A natural question is whether the reciprocal holds:
\begin{adjustwidth}{3mm}{3mm}
\textbf{Question}~(Lohrey-Steinberg)\textbf{.} Does there exist a finitely generated group with decidable Submonoid Membership and undecidable Rational Subset Membership?
\end{adjustwidth}
Lohrey and Steinberg proved that both problems are recursively equivalent in RAAGs \cite{Lohrey_Steinberg_RAAG} and infinitely-ended groups \cite{Lohrey_Steinberg_ended}. However, they conjecture a positive answer for more general groups. They also note that the existence of such a group is equivalent to the property \say{$\SMP(G)$ is decidable} not being closed under free products \cite[\S 4]{Lohrey_Steinberg_ended}.

\medbreak

For (non virtually abelian) nilpotent groups, the full picture is not clear yet.
\begin{itemize}[leftmargin=8mm]
	\item The Knapsack problem (hence Rational Subset Membership) is undecidable in large nilpotent groups, most notably $H_3(\Z)^k$ and $N_{2,k}$ (the free $2$-step nilpotent group of rank $k$) for $k\gg 1$. \cite{Knapsack_uni,Knapsack_Heis,Knapsack_Heisb}
	\item Submonoid Membership is undecidable in $H_3(\Z)^k$ for $k\gg 1$. \cite{roman2023undecidability}
\end{itemize}
Both results rely on the negative solution to Hilbert's 10th problem: there exists no algorithm deciding whether a Diophantine equation (or system of equations) admits an integer solution \cite{Hilbert10}. On the positive side, the list of results is even shorter:
\begin{itemize}[leftmargin=8mm]
	\item The Knapsack problem is decidable in $H_{2m+1}(\Z)$ for all $m\ge 1$ \cite{Knapsack_Heis}, where
	\[ H_{2m+1}(\Z) = \left\{
	\begin{pmatrix}
		1 & \mathbf a & c \\
		& I_{m} & \mathbf b^t \\
		& & 1
	\end{pmatrix} \;\middle|\; \mathbf a,\mathbf b\in \Z^m,\;c\in\Z\right\} \le SL_m(\Z). \]
	\item Colcombet, Ouaknine, Semukhin and Worrell proved that Submonoid Membership is decidable in $H_{2m+1}(\Z)$ for all $m\ge 1$. \cite{submonoid_Heis}
\end{itemize}
Note that the former relies on deep results on quadratic Diophantine equations \cite{Grunewald_Segal}, whereas the latter is elementary. This points to $\RSMP(G)$ being harder than $\SMP(G)$, hence the hope to separate both problems within the class of nilpotent groups.

\bigskip

In the first part of the paper, we re-interpret and extend Colcombet \emph{et al.}'s result:

\begin{customthm}{\ref{thm:semi}}
	There exists an algorithm with the following specifications
	\begin{itemize}[leftmargin=18mm, rightmargin=2mm]
		\item[\normalfont{Input:}\hspace*{3mm}] A finitely presented nilpotent group $G$ (given by a finite presentation), a finite set $S\subset G$ and an element $g\in G$ (given as words). \vspace*{-2mm}
		\item[\normalfont{Output:}] Finitely many instances $g_i\overset?\in R_i$ of Rational Subset Membership in a subgroup $H\le G$ such that $g\in S^*$ if and only if $g_i\in R_i$ for some $i$.
	\end{itemize}
	Moreover, it solves these instances if $h([H,H])=h([G,G])$, with $h$ the \emph{Hirsch length}.
\end{customthm}
Note that, if $h([G,G])\le 1$, e.g.\ for $H_{2m+1}(\Z)$, then either $h([H,H])=h([G,G])$ or $H$ is virtually abelian. In both cases the instances of $\RSMP(H)$ are decidable. In a more conceptual direction, Theorem \ref*{thm:semi} confirms the conjecture of Lohrey and Steinberg:
\begin{customcor}{\ref{cor:exe}}
	There exists a nilpotent group of class $2$ with decidable Submonoid Membership and undecidable Rational Subset Membership.
\end{customcor}

\begin{rem*}
	It is crucial to consider the \emph{uniform} version of the Rational Subset Membership in the previous result. Indeed the instances of $\RSMP(H)$ we need to solve \emph{will depend} on $g\in G$, even for a fixed submonoid $S^*\subseteq G$.
\end{rem*} 

\bigbreak

The second part of the paper focuses on the discrete Heisenberg group
\[ H_3(\Z) 
\simeq \la x,y \;\big|\; [x,[x,y]]=[y,[x,y]]=1 \ra. \]
The Heisenberg group is the smallest (non virtually-abelian) f.g.\ nilpotent group, in the sense that it embeds in all such groups. Our main result is the following:
\begin{customthm}{\ref{thm:rat_memb}}
	$H_3(\Z)$ has decidable Rational Subset Membership.
\end{customthm}

\medbreak

Combining this result with Theorem \ref*{thm:semi}, we get the following:
\begin{cor}
The filiform $3$-step nilpotent group {\normalfont(}also called \say{the Engel group}{\normalfont)}
	\[ E = \la x,y_1,y_2,y_3 \;\big|\; [x,y_i]=y_{i+1} \text{ for }i=1,2;\; [x,y_3]=[y_i,y_j]=1 \text{ for }i,j=1,2,3\ra \]
has decidable Submonoid Membership.
\end{cor}
Indeed, this nilpotent group has Hirsch length $h(E)=4$, hence any infinite-index subgroup admits a finite-index subgroup isomorphic to $1$, $\Z$, $\Z^2$, $\Z^3$ or $H_3(\Z)$. It can be seen as a subgroup of the unitriangular matrices $UT_4(\Z)$:
\[ E \simeq \Z^3\rtimes_X\Z = \left\{\begin{pmatrix} X^n & \mathbf y \\ 0 & 1 \end{pmatrix} \;\middle|\; n\in\Z ,\;\mathbf y\in\Z^3 \right\} \quad\text{where }X=\begin{pmatrix} 1 & 1 & 0 \\ 0 & 1 & 1 \\ 0 & 0 & 1\end{pmatrix}. \]

\newpage

We expect that some extension of Theorem \ref*{thm:rat_memb} holds:

\textbf{Conjecture A.} $H_{2m+1}(\Z)$ has decidable Rational Subset Membership for all $m\ge 1$.

In turn, this would imply that Submonoid Membership is decidable in f.g.\ nilpotent groups $G$ satisfying $h([G,G])\le 2$, or commensurable to $N_{2,3}\times\Z^n$ or $N_{3,2}\times\Z^n$. (Those groups satisfy $h([G,G])=3$, and contain no subgroup with $h([H,H])=2$.) 

\bigbreak

The proof of Theorem \ref*{thm:rat_memb} relies on a technical proposition. We reduce the problem to the class of bounded regular languages, defined by the following equivalent conditions:
\begin{thm*}[{See eg.\ \cite{Tits_for_languages}}] Fix  $\Lc\subseteq\Sigma^\star$ a regular language. The following are equivalent:
	\begin{enumerate}[leftmargin=8mm, label={\normalfont(\alph*)}]
		\item $\Lc$ has polynomial growth: $\beta_\Lc(n)\coloneqq\#\Lc \cap \Sigma^{\le n}\preceq n^d$ for some $d\ge 0$.
		\item $\Lc$ is \emph{bounded}, i.e., $\Lc \subseteq \{w_1\}^*\{w_2\}^*\ldots\{w_r\}^*$ for some $w_i\in \Sigma^\star$.
		\item $\Lc$ is a finite union of languages $t_0\,\{u_1\}^*\, t_1\, \{u_2\}^*\, t_2\ldots \{u_s\}^*\, t_s$ with $t_i,u_i\in \Sigma^\star$.
	\end{enumerate}
\end{thm*}
\smallbreak
More precisely, we prove the following result:
\begin{customprop}{\ref{thm:reduc}}
	Let $R\subseteq H_3(\Z)$ be a rational subset, i.e., $R=\ev(\Lc)$ for some regular language $\Lc\subseteq G^\star$ (see \S\ref{ssec:def} for definitions). There exists a bounded regular language $\Lc'\subset G^\star$ such that $R=\ev(\Lc')$. Moreover, $\Lc'$ can be effectively computed from $\Lc$.
\end{customprop}
This essentially reduces $\RSMP(H_3(\Z))$ to $\KP(H_3(\Z))$, which is decidable by \cite{Knapsack_Heis}.

\bigbreak

\textbf{Added in proof.} A few weeks after this article was first made public, Markus Lohrey told me about unpublished work of Doron Shafrir (2018) proving that $\SMP(G)$ is decidable and $\RSMP(G)$ is undecidable for $G=A\wr\Z^2$, with $A\ne 1$ finite abelian. This work has now been written down in Potthast's bachelor thesis \cite{Potthast}, and generalised by Ruiwen Dong in \cite{Dong_SMM}.

Moreover, building on work of Shafrir \cite{shafrir2024index}, Dong proved that the decidability of $\SMP$ is not preserved by finite-index overgroups \cite{Dong_SMM}. This leaves open the following basic questions:

\textbf{Question B.} Does there exist a finitely generated group $G$ such that
\begin{itemize}[leftmargin=6mm, label=\textbullet]
	\item $\SMP(G)$ is decidable and $\KP(G)$ is undecidable?
	\item $\RSMP(G)$ (resp.\ $\SMP$, $\KP$) is decidable and $\RSMP(G\times\Z)$ (resp.\ $\SMP$, $\KP$) is undecidable?
\end{itemize}


Finally Shafrir proved the above conjecture for rational subsets of a specific shape  \cite{shafrir2024bounded}. This is sufficient to conclude that $\SMP(G)$ is decidable for $G$ nilpotent satisfying $h([G,G])\le 2$, as well as $G=N_{2,3}\times \Z^n$ and $N_{3,2}\times\Z^n$. Extending this result further would require new ideas. However, we advertise the following problem, which still seems tractable (as the algorithms proposed in \cite{submonoid_Heis} or here do not rely on solving quadratic Diophantine equations).

\textbf{Problem C.} Give an upper bound on the time complexity of $\SMP(H_3(\Z))$.

\newpage

\section{Background on rational subsets}

\subsection{Definition} \label{ssec:def}
	Fix a group $G$. An \emph{automaton} over $G$ is a tuple $M=(V,\delta,v_0,\accept)$ with
	\begin{itemize}[leftmargin=8mm, label=-]
		\item $V$ is a finite set of vertices/states.
		\item $\delta\subset V\times G\times V$ is a finite set of edges/transitions. Each element $(u,g,v)\in\delta$ should be thought as an oriented edge $u\to v$ with label $g$.
		\item $v_0\in V$ is the initial/start vertex.
		\item $\accept\subseteq V$ is the set of terminal/accept vertices.
	\end{itemize}
	$\star$ Each automaton $M$ recognizes a language $\Lc(M)\subseteq G^\star$, namely the set of words $w\in G^\star$ we can read along (oriented) paths from $v_0$ to some accept vertex.
	
	$\star$ There is a natural map $\ev\colon G^\star\to G\colon w\mapsto \bar w$, interpreting each word as a product in $G$. 

	$\star$ A subset $R\subseteq G$ is \emph{rational} if there exists an automaton $M$ over $G$ such that $\ev(\Lc(M))=R$.
	The subset $R$ is \emph{unambiguously rational} if we can moreover ensure that $\ev\colon \Lc(M)\to R$ is bijective. In that case, $\Lc(M)$ is a \emph{regular normal form} for $R$.
	
	$\star$ An automaton is \emph{trim} if every vertex $p\in V$ lies on a path from $v_0$ to an accept vertex.
	
	$\star$ An automaton is \emph{deterministic} if, for every pair $(u,g)\in V\times G$, there exists at most one $v\in V$ such that $(u,g,v)\in \delta$. This implies that each $w\in\Lc$ is accepted by a unique path.
	
	$\star$ Given a word $w\in G^\star$, we denote its length by $\ell(w)$. Moreover, given $S\subseteq G$ and $g\in \la S\ra$, we define its \emph{word length} as $\norm{g}_S\coloneqq\min\{\ell(w)\mid w\in S^\star\text{ such that }\bar w=g\}$.

\subsection{Going to subgroups}

We state a lemma due to Gilman, himself inspired by Stalling. We recall the proof since the sets $X$ and $\tilde\Lc$ defined in the proof will turn out useful.
\begin{prop}[{\cite[Lemma 5]{Gilman1987GROUPSWA}}] \label{thm:Gil}
	Let $G$ be a group and $R\subseteq G$ a rational subset recognized by an automaton $M$ (over $G$). Suppose that $R$ sits inside a subgroup $H\le G$, then
	\begin{itemize}[leftmargin=8mm]
		\item There exists an $H$-automaton $\tilde M$ recognizing $R$ (that is, $R$ is rational in $H$). Moreover $\tilde M$ is effectively computable from $M$.
		\item If $\Lc$ is a normal form for $R$, then so is $\tilde\Lc=\Lc(\tilde M)$.
		\item If $\Lc$ is bounded, then so is $\tilde\Lc$.
	\end{itemize}
\end{prop}
In particular, from the first point, we may drop the \say{over $G$}.
\begin{proof} We start with an automaton $M=(V,\delta,v_0,\accept)$ recognizing $R$. Using the powerset construction \cite{Rabin_Scott}, we may suppose that this automaton is deterministic and trim. Let
	\begin{align*}
	X & = \left\{ \bar t\, \bar u\, \bar t^{-1}\;\middle|\; \exists p\in V \text{ s.t.} \begin{array}{c} t\in G^\star \text{ labels a simple path }v_0\to p \\ u\in G^\star \text{ labels a simple cycle }p\to p \\ \text{both paths only intersect at }p \end{array} \right\} \\[.5mm]
	Y & = \left\{ \bar y \;\big|\; \exists q\in\accept \text{ s.t. }y\in G^\star \text{ labels a simple path }v_0\to q \right\}
	\end{align*}
	We consider the regular language
	\[ \tilde \Lc = \left\{ x_1x_2\ldots x_\ell\, y \in X^*Y \;\middle|\; \begin{array}{l} t_i\text{ is a prefix of }t_{i+1}u_{i+1} \text{ for }1\le i\le \ell-1 \\ t_\ell\text{ is a prefix of }y,\;\text{ where }\bar t_i\bar u_i\bar t_i^{-1}=x_i \end{array} \right\} \] 
	which is recognized by the automaton $\tilde M=(\tilde V,\tilde\delta, \varepsilon,\{\anode\})$ with
	\begin{itemize}[leftmargin=8mm, label=-]
		\item $\tilde V= \{t\in G^\star \mid \exists p\in V \text{ s.t. }t\text{ labels a simple path }v_0\to p\} \sqcup \{\anode\}$,
		\item $\tilde\delta$ consists of all edges $s\to t$ labeled by $x\in X$ if $x=\bar t\,\bar u\,\bar t^{-1}$ and $s$ is a prefix of $tu$, \newline
		and all an edge $t\to \anode$ labeled by $\bar y\in Y$ if $t$ is a prefix of $y$,
		\item the empty string $\varepsilon$ (which labels the simple path $v_0\to v_0$) is the starting vertex,
		
		\item and the vertex $\anode$ is the unique accept vertex.
	\end{itemize}
	
	\smallbreak
	
	$\blacktriangleright$ We first prove that $\ev(\Lc)=\ev(\tilde\Lc)$.
	
	Any word $w\in \Lc$ is accepted by a unique path $v_0\to v_0$ in the automaton $M$. This path can be decomposed as a product of conjugates $t_iu_it_i^{-1}$, using a loop-erasure algorithm:
	\begin{center}
		\begin{tikzpicture}[scale=1.05]
	\draw[black, thick, -latex] (-1,-1)
		to [out=60, in=-180, looseness=.9] (0,0)
		to [out=0, in=20, looseness=2] (-.5,1)
		to [out=-160, in=-160, looseness=3] (-.8,1.5)
		to [out=20, in=100, looseness=2] (-.5,1)
		to [out=-80, in=120] (0,0)
		to [out=-60, in=170] (0.5,-.5)
		to [out=-10, in=-85, looseness=2] (1,-.2)
		to [out=95, in=90, looseness=2] (0.5,-.5)
		to [out=-90, in=0, looseness=1.5] (-.4,-1);
	\node at (-.8,-.1) {$w$};
	\node at (1.4,0.3) {$=$};
	
	\draw[TealBlue, -latex, thick, shift={(3,0)}] (-1,-1)
		to [out=60, in=-180, looseness=.9] (0,0)
		to [out=0, in=20, looseness=2] (-.5,1);
	\draw[black, thick, -latex, shift={(3,0)}] (-.5,1)
		to [out=-160, in=-160, looseness=3] (-.8,1.5)
		to [out=20, in=100, looseness=2] (-.5,1);
	\node[TealBlue] at (2.2,-.1) {$t_1$};
	\node at (2.8,1.5) {$u_1$};
	\node at (3.7,0.3) {$\cdot$};
	
	\draw[TealBlue, -latex, thick, shift={(5,0)}] (-1,-1)
		to [out=60, in=-180, looseness=.9] (0,0);
	\draw[black, thick, -latex, shift={(5,0)}] (0,0)
		to [out=0, in=20, looseness=2] (-.5,1)
		to [out=-80, in=120] (0,0);
	\node[TealBlue] at (4.2,-.1) {$t_2$};
	\node at (5.1,1.2) {$u_2$};
	\node at (5.8,0.3) {$\cdot$};
	
	\draw[TealBlue, -latex, thick, shift={(7,0)}] (-1,-1)
		to [out=60, in=-180, looseness=.9] (0,0)
		to [out=-60, in=170] (0.5,-.5);
	\draw[black, thick, -latex, shift={(7,0)}] (0.5,-.5)
		to [out=-10, in=-85, looseness=2] (1,-.2)
		to [out=95, in=90, looseness=2] (0.5,-.5);
	\node[TealBlue] at (6.2,-.1) {$t_3$};
	\node at (7.5,0.2) {$u_3$};
	\node at (8.4,0.3) {$\cdot$};
	
	\draw[black, thick, -latex, shift={(9.7,0)}] (-1,-1)
		to [out=60, in=-180, looseness=.9] (0,0)
		to [out=-60, in=170] (0.5,-.5)
		to [out=-90, in=0, looseness=1.5] (-.4,-1);
	\node at (9.7,.2) {$y$};
\end{tikzpicture}

		\captionsetup{margin=10mm, font=footnotesize}
		
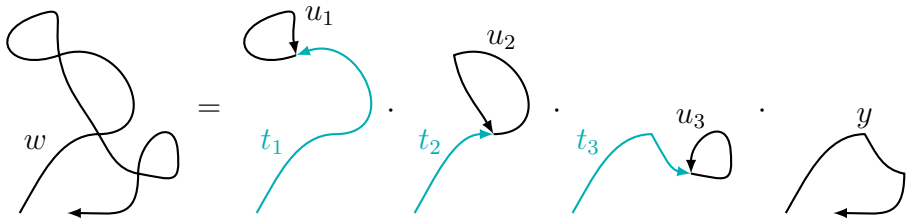
\captionof{figure}{A path decomposed as a product of \say{freeze frames} of the loop-erasure algorithm.}
	\end{center}
	
	The resulting word $\tilde w=(\bar t_1\bar u_1\bar t_1^{-1})\ldots (\bar t_\ell\bar u_\ell\bar t_\ell^{-1})\cdot \bar y\in \tilde\Lc$ is the \emph{decomposition} of $w$. Both $w$ and its decomposition $\tilde w$ evaluate to the same element in $G$. Moreover, words in $\tilde\Lc$ are exactly decompositions of words in $\Lc$. This proves that $\ev(\Lc)=\ev(\tilde\Lc)$.
	
	\medbreak
	
	$\blacktriangleright$ We prove that $X,Y\subseteq H$.
	
	Take $x=\bar t\,\bar u\,\bar t^{-1}\in X$, where $t$ labels a simple path $v_0\to p$ and $u$ labels a simple cycle $p\to p$ with $t$ and $u$ only intersecting at $p$. As the automaton $M$ is trim, there exists a path $p\to q\in \accept$ labeled by some word $v\in G^\star$. Consider $tuv\in\Lc$ and $xx_2\ldots x_\ell y\in\tilde\Lc$ its decomposition. Observe that $x_2\ldots x_\ell y\in\tilde\Lc$ too, so
	\[ x = (xx_2\ldots x_\ell y)(x_2\ldots x_\ell y)^{-1} \in R\cdot R^{-1} \subseteq H. \]
	On the other side, we have $Y\subseteq R\subseteq H$.
	
	\medbreak
	
	$\blacktriangleright$ If $\Lc$ is a normal form, then each $g\in R$ is represented by a unique word $w$ in $\Lc$, hence has a unique representative $\tilde w$ in $\tilde\Lc$.
	
	\medbreak
	
	$\blacktriangleright$ Observe that $\ell(w)\le \abs V\ell(\tilde w)$ for all $w\in\Lc$ (as each letter in $\tilde w$ corresponds to a simple cycle or a simple path in $w$), hence $\beta_\Lc(n)\succeq\beta_{\tilde\Lc}(n)$. If $\Lc$ is bounded, then $\Lc$ and $\tilde\Lc$ have polynomial growth, hence $\tilde\Lc$ is bounded \cite{Tits_for_languages}.
\end{proof}

\newpage

Another lemma which complements nicely is the following:

\begin{lem} \label{lem:FinInd_inter} Let $G$ be a group and $H\le G$ is a finite-index subgroup.
	\begin{enumerate}[leftmargin=8mm, label={\normalfont(\alph*)}]
		\item  If $R\subseteq G$ is rational, then $Hg\cap R$ is rational for each $g\in G$.
	\end{enumerate}
	If moreover we are given an automaton $M$ such that $R=\ev(\Lc(M))$ with labels inside a finite set $S\subseteq G$, and the Schreier graph $\Sch(H\backslash G,S)$, then 
	\begin{enumerate}[leftmargin=8mm, label={\normalfont(\alph*)}]
		\setcounter{enumi}{1}
		\item we can effectively compute an automaton for $Hg\cap R$, and
		\item we can effectively compute $\left\{Hg\in H\backslash G \;\big|\; Hg\cap R\neq\emptyset\right\}$ and compute a subset $\Kc\subset \Lc(M)$ consisting of a unique representative for each of these cosets. 
	\end{enumerate}
\end{lem}
\begin{proof}
	(a-b) Given an automaton for $M=(V,\delta,v_0,\accept)$ for $R$, we construct the following automaton recognizing $Hg\cap R$:
	\begin{itemize}[label=-, leftmargin=8mm]
		\item $V'=H\backslash G\times V$
		\item $v_0'=(H,v_0)$ 
		\item $\accept'=\{Hg\}\times\accept$
		\item $(Hg_1,v_1)\overset s\to (Hg_2,v_2)$ if and only if $Hg_1s=Hg_2$ and $v_1\overset s\to v_2$.
	\end{itemize}
	
	\medbreak
	
	(c) For each coset, one can compute an automaton for $Hg\cap R$ and then decide if the intersection is non-empty: does there exist a path $v_0'\to \accept'$? If yes, a representative is given by $w\in \Lc(M)$ labeling the first path found.
\end{proof}
\begin{rem}
	Whenever $G=\la A\mid B\ra$ is finitely presented, the Todd-Coxeter algorithm allows to compute $\Sch(H\backslash G,S)$ taking as input $S$ and generators for $H$ (given as words over $A$), so $\Sch(H\backslash G,S)$ doesn't need to be part of the input of our algorithm.
\end{rem}

\subsection{Virtually abelian groups}

We recall a result of Eilenberg and Schützenberger.
\begin{thm}[{\cite{EILENBERG1969}}] \label{sec1:Eilenberg}
	Let $G$ be an abelian group and $R$ be a rational subset, then $R$ is unambiguously rational. Moreover, given a language $\Lc$ such that $\ev(\Lc)=R$, we can effectively compute a regular normal form $\Lc'$ for $R$.
\end{thm}
This result can easily be extended to virtually abelian groups using Lemma \ref{lem:FinInd_inter}(b).
\begin{cor} \label{cor:bounded_in_virtab}
	Any rational subset of a virtually abelian group can be represented by a bounded regular language $\Lc'$. Moreover, $\Lc'$ can be computed effectively.
\end{cor}
\begin{proof}
	Take $\Lc'$ the normal form from Theorem \ref*{sec1:Eilenberg}, and let $S$ be the (finite) set of elements appearing as labels on an automaton recognizing $\Lc'$. Note that $\ell(w)\ge \norm{\bar w}_S$ for all $w\in S^*$ hence we can compare the volume growth of $\Lc'$ and $(\la S\ra,S)$: we have $\beta_{\Lc'}(n)\le\beta_{(\la S\ra,S)}(n)$, but $\la S\ra$ is a finitely generated virtually abelian group, hence has polynomial growth, so the rational language $\Lc'$ has polynomial growth hence is bounded (see \cite{Tits_for_languages}).
\end{proof}
\section{Dimension gain for the Submonoid Membership}

\subsection{Background on nilpotent groups}

We first recall the classical notion of Hirsch length, and some related properties.

\begin{defi} Let $G$ be a group. The \emph{lower central series} is the sequence of subgroups
	\[ \gamma_1(G)=G,\quad \gamma_{i+1}(G) = [\gamma_i(G),G] \;\,\text{for all }i\ge 1 \]
	We say that $G$ is \emph{nilpotent of class $c$} (or $c$-step nilpotent) if $\gamma_{c+1}(G)=\{e\}$.
\end{defi}
\begin{defi} The \emph{Hirsch length} of a $c$-step nilpotent group $G$ is defined as
	\[ h(G) = \sum_{i=1}^c \rk_\Q\!\big( \gamma_i(G)/\gamma_{i+1}(G)\big). \]
\end{defi}
\begin{prop}[{Folklore, \cite[Exercise 8]{Segal_1983}}] \label{prop:Hirsh}
	Let $G$ be a finitely generated nilpotent group.
\begin{enumerate}[label={\normalfont(\alph*)}, leftmargin=8mm]
	\item For any exact sequence $1\to N\to G\to Q\to 1$, we have $h(G)=h(N)+h(Q)$.
	\item For any group $H\le G$, we have $h(H)\le h(G)$ with equality if and only if $[G:H]<\infty$.
	\item Let $N_{2,m}$ be the free $2$-step nilpotent group on $m$ generators. Then $h(N_{2,m})=m+{m\choose 2}$.
\end{enumerate}
\end{prop}

\bigbreak

Next we quote a useful result, which can be understood as a discrete version of the Chow-Rashevskii theorem in sub-Riemannian geometry. This requires a definition/lemma:
\begin{lem}[Torsionfree part of the abelianization] \label{lem:tf_abelian} Let $G$ be a group, then
	\[ \overline{[G,G]}= \left\{g\in G \;\big|\; \exists n\ne 0,\; g^n \in [G,G]\right\} \]
	is a characteristic subgroup of $G$. Moreover, if $G$ is finitely generated, then $[G,G]$ has finite index inside $\overline{[G,G]}$, and $G/\overline{[G,G]}\simeq\Z^r$ for some $r$.
\end{lem}
\begin{proof}
	Let $\tau\colon G\to G/[G,G]$ be the abelianization map, then torsion elements of $G/[G,G]$ form a characteristic subgroup $T$ and $\overline{[G,G]}=\tau^{-1}(T)$. If $G$ is finitely generated, then $G/[G,G]\simeq \Z^r\times T$ hence $[\overline{[G,G]}:[G,G]]=\abs T$ and $G/\overline{[G,G]}\simeq\Z^r$.
\end{proof}
\begin{prop}[{\cite{BCM24}, see also \cite{shafrir2024saturation}}] \label{prop:BCM}
	Consider an infinite, finitely generated nilpotent group $G$, and $\pi\colon G\onto G/\overline{[G,G]}\simeq\Z^r$. For any $S\subseteq G$, the following are equivalent
	\begin{enumerate}[leftmargin=8mm, label={\normalfont(\alph*)}]
		\item $\Conv(\pi(S))\subseteq \R^r$ contains a ball $B(\mathbf 0,\varepsilon)$ for some $\varepsilon>0$.
		\item For every non-zero homomorphism $f\colon G\to \R$, there exists $s\in S$ s.t.\ $f(s)<0$.
		\item The submonoid $S^*$ is a finite-index subgroup of $G$.
	\end{enumerate}
	If $S$ is finite, we can restrict to non-zero homomorphism $f\colon G\to\Z$.
\end{prop}

\medbreak

\subsection{Proof of Theorem \ref*{thm:semi}}

We prove the following algorithmic reduction:
\begin{thm}\label{thm:semi}
There exists an algorithm with the following specifications:
\begin{itemize}[leftmargin=18mm, rightmargin=2mm]
	\item[\normalfont{Input:}\hspace*{3mm}] A finitely presented nilpotent group $G$, a finite set $S\subset G$, and $g\in G$. 
	\item[\normalfont{Output:}] Finitely many instances $\{g_i\overset?\in R_i\}$ of the Rational Subset Membership in a subgroup $H\le G$ such that $g\in S^*$ if and only if $g_i\in R_i$ for some $i$.
\end{itemize}
Moreover, the algorithm solves these instances if $h([H,H])=h([G,G])$.
\end{thm}

\begin{proof} First, we may assume $e_G\in S$. The proof splits into three steps
	\begin{enumerate}[leftmargin=8mm, label=(\arabic*)]
		\item We compute the image of $S$ through a map $\pi\colon G\onto\Z^r$ as in Proposition \ref{prop:BCM}. Using this quotient, we define a subgroup  $H\le G$ (depending only on $S$) and a partition $S=S_0\sqcup S_+$, looking whether $\pi(s)$ is invertible or not in $\pi(S)^*\le\Z^r$ for each $s\in S$.
		\item We reduce the problem $g\in S^*$ to finitely many instances of $\RSMP(H)$. 
		\item In the case when $h([H,H])=h([G,G])$, we solve the previous instances of $\RSMP(H)$. 
	\end{enumerate}
	The first step follows \cite[Theorem 7]{submonoid_Heis} closely, while the last step generalizes \say{Case II} from the same proof. The observation of the second step seems new.

\bigskip

$\blacktriangleright$ First we compute a map $\pi\colon G\to \Z^r$ from a presentation for $G$. We consider the polytope $\Conv(\pi(S))\ni\mathbf 0$ given by a V-representation (namely, $\pi(S)$). We can compute a H-representation (H for half-space), that is, a finite set of inequalities $\{f_i(\mathbf v)\ge a_i \mid i\in I \}$ with $f_i\colon\Z^r\to\Z$ non-zero linear forms and $a_i\in\Z$ such that
\[
\Conv(\pi(S)) = \left\{ \mathbf v\in\R^r \;\big|\; \forall i\in I,\; f_i(\mathbf v)\ge a_i \right\}.
\]
This is the classical \emph{facet enumeration problem}. (See \cite[Section 1.2]{Ziegler} and references therein.) We compute a maximal linearly independent set $\{f_1,f_2,\ldots,f_s\}$ inside $\{f_i\mid i\in I, a_i=0\}$, and define a homomorphism $f\colon G\to\Z^s$ via
\[ f(h) = \big(f_1(\pi(h)),\,f_2(\pi(h)),\ldots,\,f_s(\pi(h))\big). \]
The image of $f\colon G\to\Z^s$ has finite index (by linear independence), and $f(S)\subset \Z_{\ge 0}^s$. We define $H=\ker f$ and partition $S$ into $S_0=S\cap H$ and $S_+=S\setminus H$.

\begin{adjustwidth}{3mm}{3mm}
\textbf{Geometric intermezzo:} We are looking at the minimal face $F$ of $\Conv(\pi(S))$ containing $\mathbf 0$ (i.e., no proper sub-face of $F$ contains $\mathbf 0$). This face is given by
\[ F = \Conv(\pi(S)) \cap K=\Conv(\pi(S_0))\]
where $K=\left\{\mathbf v\in\R^r \;\big|\; f_i(\mathbf v)=0 \text{ for }i=1,2,\ldots,s\right\}$.
The minimality of $F$ can be stated as \say{there exists $\varepsilon>0$ such that $B(\mathbf 0,\varepsilon)\cap K\subset F$}.
\end{adjustwidth}
\medbreak

$\blacktriangleright$ There exists only finitely many words $w=u_1u_2\ldots u_k\in S_+^\star$ such that 
	\[ f(\bar w) = \sum_{i=1}^k f(u_i) = f(g) \]
	(Indeed, $k$ is bounded by the sum of the components of $f(g)$.) For each word $w$, we need to decide whether $g \in S_0^* u_1 S_0^* u_2S_0^* \ldots S_0^* u_k S_0^*$. Observe that
	\[ S_0^* \cdot u_1 S_0^* \cdot u_2S_0^* \cdot \ldots \cdot u_k S_0^* = S_0^* \cdot \big(v_1S_0v_1^{-1}\big)^* \cdot  \big(v_2S_0v_2^{-1}\big)^* \cdot \ldots \cdot \big(v_kS_0v_k^{-1}\big)^* \cdot v_k \]
	where $v_i=u_1u_2\ldots u_i$. Hence the problem can be restated as
	\[ gv_k^{-1} \overset?\in R \coloneqq  S_0^*  \cdot \big(v_1S_0v_1^{-1}\big)^* \cdot \big(v_2S_0v_2^{-1}\big)^* \cdot  \ldots\cdot  \big(v_kS_0v_k^{-1}\big)^* \subseteq H \]
	which is an instance of the Rational Subset Membership Problem in $H$.

	We note that algorithms presented in \cite{subgroup_pres} allow to compute a presentation for $H$, and then rewrite all elements $gv_k^{-1}\in H$ and $v_isv_i^{-1}\in H$ as words over the corresponding generating set. (Problems (III), (IV) and (II) in the article.)

\bigbreak

$\blacktriangleright$ We solve the Membership problem under the hypothesis $h([G,G])=h([H,H])$.

Observe that $\overline{[G,G]}=\overline{[H,H]}$. Indeed, we have $\overline{[G,G]}\le H$ as $H$ is the kernel of $f\colon G\to\Z^s$ with $\Z^s$ torsionfree abelian. Moreover Proposition \ref{prop:Hirsh}(b) and Lemma \ref{lem:tf_abelian} imply that
\[ m = \Big[ \overline{[G,G]} : [H,H] \Big] = \Big[ \overline{[G,G]} : [G,G] \Big] \Big[[G,G] : [H,H] \Big] < \infty. \]
For every $g\in \overline{[G,G]}$, we have $g^{m!}\in[H,H]$ hence $g\in\overline{[H,H]}$, proving the observation. In particular, we can identify the quotient $H/\overline{[H,H]}=H/\overline{[G,G]}=\pi(H)$.

\medbreak

Condition (a) of Proposition \ref{prop:BCM} (for $H$ and $S_0$) now reads \say{$F=\Conv(\pi(S_0))$ contains a ball $B(\mathbf 0,\varepsilon)\cap K$}, which holds by minimality of $F$. We conclude that $S_0^*$ is a finite-index subgroup of $H$. Finally, the problem can be restated as
	\[ \la S_0\ra gv_k^{-1} \;\cap\; \big(v_1S_0v_1^{-1}\big)^* \cdot \big(v_2S_0v_2^{-1}\big)^* \cdot  \ldots\cdot  \big(v_kS_0v_k^{-1}\big)^* \overset? =\emptyset \]
which is easily decided using Lemma \ref{lem:FinInd_inter}(c).	
\end{proof}

\subsection{Proof of Theorem \ref*{cor:exe}}

We give a lemma on abstract commensurability classes of subgroups of $N_{2,m}\times\Z^n$. Some ideas can already be found in \cite[Theorem 7]{golovin1955subgroups}.
\begin{lem}\label{lem:subgroup_of_free}
	Any subgroup $H\le N_{2,m}\times\Z^n$ admits a finite-index subgroup $H'\trianglelefteq H$ isomorphic to $N_{2,k}\times\Z^\ell$ for some $k\le m$ and $\ell\le {m\choose 2}-{k\choose 2}+n$. 
\end{lem}
\begin{proof}
	
	Let $G=N_{2,m}\times\Z^n$ and fix a subgroup $H\le G$. Consider the composition $\alpha=\alpha_2\circ \alpha_1$ where $\alpha_1$ is the projection on the first factor, and $\alpha_2$ is the abelianization map:
	\begin{center}
		\begin{tikzcd}
			N_{2,m}\times\Z^n \arrow[r, two heads, "\alpha_1"] 
			& N_{2,m} \arrow[r, two heads, "\alpha_2"]
			& \Z^m
		\end{tikzcd}
	\end{center}
	The image $\alpha(H)\le\Z^m$ is isomorphic to $\Z^k$ for some $k\le m$. Fix $g_1,\ldots,g_k\in H$ such that $\{\alpha(g_1),\ldots,\alpha(g_k)\}$ is a basis for $\alpha(H)$.
	\begin{adjustwidth}{3mm}{3mm}
		\textbf{Claim 1.} $\la g_1,\ldots,g_k\ra\simeq N_{2,k}$ and $\{g_1,\ldots,g_k\}$ is a basis.
		
		Observe that the projections of $\alpha_1(g_1),\ldots,\alpha_1(g_k)\in N_{2,m}$ in $\Z^m$ are linearly independent, hence they generate a subgroup isomorphic to $N_{2,k}$ as a basis (see \cite[Theorem 1.3]{Moran_Sub_of_Free}). Moreover, since $N_{2,k}$ is relatively free in the variety of $2$-step nilpotent groups and $N_{2,m}\times\Z^n$ is $2$-step nilpotent, we can define a morphism $\alpha_1(g_i)\mapsto g_i$, proving that the restriction $\alpha_1|_{\la g_1,\ldots,g_k\ra}$ is invertible, hence the claim.
	\end{adjustwidth}
	
	Observe that $K=\ker \alpha$ is the center of $G$, in particular is abelian. We have
	\begin{center}
		\begin{tikzcd}
			K \arrow[r, two heads] 
			& K / \big(\!\la g_1,\ldots,g_k\ra\cap K\big) \arrow[r, "\sim"]
			& T\times \Z^{{m\choose 2}-{k\choose 2}+n} \arrow[r, two heads]
			& \Z^{{m\choose 2}-{k\choose 2}+n}
		\end{tikzcd}
	\end{center}
	where $T$ is some finite abelian group.  Indeed $K\simeq \Z^{{m\choose 2}+n}$ and $\la g_1,\ldots,g_k\ra\cap K\simeq\Z^{k\choose 2}$. We denote the whole composition by $\beta$, and the composition of the first two arrows by $\gamma$.
	
	\bigskip
	
	The image $\beta(H\cap K)\le \Z^{{m\choose 2}-{k\choose 2}+n}$ is isomorphic to $\Z^\ell$ for some $\ell\le {m\choose 2}-{k\choose 2}+n$. Fix $h_1,\ldots,h_\ell\in H\cap K$ such that $\{\beta(h_1),\ldots,\beta(h_\ell)\}$ is a basis of $\beta(H\cap K)$.
	
	\begin{adjustwidth}{3mm}{3mm}
		\textbf{Claim 2.} Since $K$ is abelian, $\beta\colon\la h_1,\ldots,h_\ell\ra\onto \beta(H\cap K)\simeq \Z^\ell$ is an isomorphism.
		
		\medbreak
		
		\textbf{Claim 3.} We have $H'\coloneqq\la g_1,\ldots,g_k,h_1,\ldots,h_\ell\ra\simeq N_{2,k}\times\Z^\ell$. It suffices to check that
		\begin{itemize}[leftmargin=6mm]
			\item $\la g_1,\ldots,g_k\ra$ and $\la h_1,\ldots,h_\ell\ra$ commutes, which is true as $\la h_1,\ldots,h_\ell\ra\le K=\Zc(G)$.
			\item $\la g_1,\ldots,g_k\ra\cap\la h_1,\ldots,h_\ell\ra=\{e\}$, which is true as restriction $\beta|_{\la g_1,\ldots,g_k\ra\cap K}$ is zero while the restriction $\beta|_{\la h_1,\ldots,h_\ell\ra}$ is injective.
		\end{itemize}
		
		\medbreak
		
		\textbf{Claim 4.} $H'\trianglelefteq H$ and $H/H'$ is finite.
		
		Recall that $K=\ker\alpha$ is central hence
		\[
		\forall g,g',h,h'\in G,\; \big(\alpha(g)=\alpha(g')\text{ and }\alpha(h)=\alpha(h')\big)\implies [g,h]=[g',h'].
		\]
		Therefore $\alpha(H)=\alpha(H')$ implies $[H,H]=[H',H']$. It follows that $H'\trianglelefteq H$ (as any subgroup $H'\le H$ containing $[H,H]$ is normal in $H$). Moreover,
		\[
		H/H'\simeq \big(H\cap K\big)/\big(H'\cap K\big) \simeq \gamma(H\cap K)/\gamma(H'\cap K) \simeq \gamma(H\cap K)\cap T
		\]
		using the Nine lemma, the third isomorphism theorem, and the Nine lemma again.
	\end{adjustwidth}
		\adjustbox{scale=.6,center}{%
			\begin{tikzcd}
			& 1 \arrow[d] & 1\arrow[d] & 1\arrow[d] & \\
			1\arrow[r] & H'\cap K \arrow[r]\arrow[d] & H'\arrow[r]\arrow[d] & \alpha(H')\arrow[r]\arrow[d, "\wr"] & 1 \\
			1\arrow[r] & H \cap K\arrow[r]\arrow[d] & H\arrow[r]\arrow[d] & \alpha(H)\arrow[r]\arrow[d] & 1 \\
			1\arrow[r] & (H\cap K)/(H'\cap K) \arrow[r,dashed,"\sim"]\arrow[d] & H/H'\arrow[r]\arrow[d] & 1 \\ & 1 & 1 & & 
			\end{tikzcd}
			\quad
			\begin{tikzcd}
				& & 1\arrow[d] & 1\arrow[d] & \\
				& 1 \arrow[r]\arrow[d] & \gamma(H'\cap K)\arrow[r, "\sim"]\arrow[d] & \beta(H'\cap K)\arrow[r]\arrow[d, "\wr"] & 1 \\
				1\arrow[r] & \gamma(H \cap K)\cap T\arrow[r]\arrow[d] & \gamma(H\cap K) \arrow[r]\arrow[d] & \beta(H\cap K)\arrow[r]\arrow[d] & 1 \\
				1\arrow[r] & \gamma(H \cap K)\cap T \arrow[r,dashed,"\sim"]\arrow[d] & \gamma(H\cap K)/\gamma(H'\cap K) \arrow[r]\arrow[d] & 1 \\ & 1 & 1 & & 
			\end{tikzcd}}
	
	It should be noted that everything is effectively computable from generators for $H$. 
\end{proof}

\newpage

Finally, we deduce the existence of a group separating $\SMP(G)$ and $\RSMP(G)$:

\begin{cor} \label{cor:exe}
	There exist $m,n\ge 0$ such that $N_{2,m}\times\Z^n$ has decidable Submonoid Membership Problem and undecidable Rational Subset Membership Problem.
\end{cor}
\begin{proof}
	Consider a group $G=N_{2,m}\times\Z^n$ with undecidable $\RSMP(G)$ and $m$ minimal. Such a group does exist by \cite{Knapsack_Heis}, and \cite{Knapsack_Heisb} even proves that $m\le 26$.
	
	We prove that $\SMP(G)$ is decidable. Using Theorem \ref{thm:semi}, any instance of $\SMP(G)$ reduces to finitely many instances of $\RSMP(H)$ for a subgroup $H\le G$ with $h([H,H])<h([G,G])$. However Lemma \ref{lem:subgroup_of_free} tells us such subgroups admit finite-index subgroups $H'\simeq N_{2,k}\times\Z^\ell$ with $k<m$ and $\ell\le {m\choose 2}-{k\choose 2}+n$. In turn, $\RSMP(H)$ reduces to $\RSMP(N_{2,k}\times \Z^\ell)$ using Lemma \ref{lem:FinInd_inter}(a-b) and Proposition \ref{thm:Gil} (see also \cite{Grunschlag}), and the latter is decidable as $k<m$.
\end{proof}
\begin{rem}
	It is important that subgroups of $N_{2,m}\times\Z^n$ fall into finitely many abstract commensurability classes, otherwise $\RSMP(H)$ could be decidable for each group $H$ without any uniform algorithm working for all $H$. Here, we can compute the values of $k$ and $\ell$, and decide which of the finitely many algorithms to apply.
\end{rem}
\section{Reduction to bounded regular languages}

\subsection{The Heisenberg group} \label{sec:Heis_model}

We recall the exponential coordinates on $H_3(\Z)$, and a geometrical interpretation for them. The group $H_3(\Z)$ is generated by the matrices
\[
x = \begin{pmatrix} 1 & 1 & 0 \\ \color{gray}{0} & 1 & 0 \\ \color{gray}{0} & \color{gray}{0} & 1 \end{pmatrix}
\quad \text{and} \quad 
y = \begin{pmatrix} 1 & 0 & 0 \\  \color{gray}{0} & 1 & 1 \\ \color{gray}{0} & \color{gray}{0} & 1 \end{pmatrix}.
\]
Another notable element is $z\coloneqq[x,y]=xyx^{-1}y^{-1}$. For a general element
\[ g = \begin{pmatrix} 1 & a & c \\ \color{gray}{0} & 1 & b \\ \color{gray}{0} & \color{gray}{0} & 1 \end{pmatrix} \in H_3(\Z), \]
we define $\hat g=(a,b)\in\Z^2$ and $A(g)=c-\frac12 ab$. Note that $\pi\colon g\mapsto \hat g$ is the abelianization map. For instance $\hat z=\0$ (since $z\in [G,G]$) and $A(z)=1$.


Given a product $g=s_1s_2\ldots s_n$ with $s_i\in\{x^\pm,y^\pm\}$, the coordinates $\hat g$ and $A(g)$ can be interpreted as follows: we consider the path $\gamma_g$ in $\Z^2$, starting from $\mathbf 0$, where the $i$-th step goes in left, right, up or down direction depending if $s_i=x$, $x^{-1}$, $y$ or $y^{-1}$. Then \vspace*{1mm}

\begin{minipage}{.68\linewidth}
\begin{itemize}[leftmargin=8mm]
	\item $\hat g$ is the endpoint of this path, and \vspace*{3mm}
	\item $A(g)$ is the algebraic area (or \emph{balayage} area) of the curve obtained by closing the path by a straight segment from $\hat g$ to $\mathbf 0$. Each point of the plane is weighted according to the winding number of the path around it. \vspace*{3.5mm} 
\end{itemize}
The adjacent figure is the path representing the product $g=y^{-1}x^3yx^{-1}y^{-2}x^3y^3x^{-3}y^3x^2$, with exponential coordinates given by $\hat g=(4,4)$ and $A(g)=12$.
\end{minipage} \hspace*{5mm}
\begin{minipage}{.3\linewidth}
	\centering
	\begin{tikzpicture}[scale=.55]
		\fill[LimeGreen!25, rounded corners=2pt] (0,0) -- (0,-1) -- (2,-1) -- (2,-2) -- (5,-2) -- (5,1) -- (2,1) -- (2,2);
		\fill[LimeGreen!50, rounded corners=2pt] (2,-1) -- (3,-1) -- (3,0) -- (2,0) -- (2,-1);
		\fill[Magenta!25, rounded corners=2pt] (2,2) -- (2,4) -- (4,4);
		
		\draw[black!50] (-1.2,-2.2) grid (6.2,5.2);
		\draw[very thick, black, -latex] (0,-2.2) -- (0,5.3);
		\draw[very thick, black, -latex] (-1.2,0) -- (6.3,0);
		
		\draw[ultra thick, Turquoise, rounded corners=2pt, -latex] (0,0) -- (0,-1) -- (3,-1) -- (3,0) -- (2,0) -- (2,-2) -- (5,-2) -- (5,1) -- (2,1) -- (2,4) -- (4,4);
		\draw[very thick, Turquoise, dashed] (0,0) -- (4,4);
		
		\node at (2.5,-.5) {\footnotesize$+2$};
		\node at (2.5,0.5) {\footnotesize$+1$};
		\node at (2.5,1.5) {\footnotesize$0$};
		\node at (2.5,3.5) {\footnotesize$-1$};
	\end{tikzpicture}
\end{minipage}

\begin{lem}[{See \cite[Proposition 4.1]{Intermediate}}]
	The exponential coordinates satisfy
	\[ \widehat{gh} = \hat g+\hat h \qquad\text{and}\qquad	A(gh) = A(g) + A(h) + \frac12 \det(\hat g;\hat h).  \]
	As a corollary, we also have
	\[ A(g_1g_2\ldots g_n) = \frac12\sum_{i<j} \det(\hat g_i;\hat g_j) + \sum_{i=1}^n A(g_i)
	\qquad\text{and}\qquad
	[g,h] = z^{\det(\hat g;\hat h)}. \]
\end{lem}
\begin{lem} \label{lem:compare_mod_N}
	Let $N\trianglelefteq G=H_3(\Z)$ be a finite-index normal subgroup. Fix $a,b\in N$ such that $\pi(N)=\langle \hat a,\hat b\rangle$ and $d\in\Z_{>0}$ such that $N\cap[G,G]=\langle z^d\rangle$. Two elements $g,h\in H_3(\Z)$ lie in the same $N$-coset if and only if the following conditions hold:
	\begin{enumerate}[leftmargin=8mm, label=(\arabic*)]
		\item $\hat g-\hat h \in \langle \hat a,\hat b\rangle$, say $\hat g-\hat h=m\hat a+n\hat b$ with $m,n\in\Z$.
		\item $A(g) \equiv A(h\cdot a^mb^n) = A(h) + A(a^mb^n) + \frac12\det(\hat h;m\hat a+n\hat b) \pmod d$.
	\end{enumerate}
	In particular, if $g=h$ in $G/(N\cap[G,G])$ (i.e., $m=n=0$), then $A(g)-A(h)\in d\Z$. For $N=G$, this means that $A(g)-A(h)\in\Z$ as soon as $\hat g=\hat h$.
\end{lem}
\begin{proof}
	Condition (1) is clearly necessary. Supposing (1), we have $a^mb^n \in N$ and $g[G,G]=ha^mb^n[G,G]$, hence $gN=hN$ is equivalent to $g\langle z^d\rangle=ha^mb^n\langle z^d\rangle$, which is condition (2).
\end{proof}
\vspace*{-2mm}

\subsection{Reduction to rational submonoids} \label{ssec:reduc}
The goal of \S\ref{ssec:reduc}, \S\ref{ssec:cases} is to prove the following:
\begin{prop} \label{propE} \label{thm:reduc}
	Let $R\subseteq H_3(\Z)$ be a rational subset, i.e., $R=\ev(\Lc)$ for some regular language $\Lc\subseteq G^\star$. There exists a \emph{bounded} regular language $\Lc'\subset G^\star$ such that $R=\ev(\Lc')$. Moreover, an automaton for $\Lc'$ can be effectively computed from an automaton for $\Lc$.
\end{prop}
We begin by reducing to the case $\start=\accept=\{v\}$ (i.e., to the case of rational submonoids). For this step, $G$ can still be an arbitrary group. We consider a general rational subset $R=\ev(\Lc)$, where $\Lc$ is accepted by a automaton $M=(V,\delta,v_0,\accept)$ over $G$.

$\blacktriangleright$ Using a loop-erasure algorithm, we decompose each word $w\in \Lc$ as
\[ w = w_0s_1w_1s_2\ldots s_\ell w_\ell \]
where $s_1,s_2,\ldots, s_\ell\in G$ label a simple path $v_0\overset{s_1}\to v_1\overset{s_2}\to \ldots \overset{s_\ell}\to v_\ell$ with $v_\ell\in\accept$, and each $w_i\in G^\star$ labels a cycle $v_i\to v_i$. We proceed as follows: start at the vertex $v_0$ and skip directly to the last visit of $v_0$, hence bypassing a (possibly empty) cycle $w_0$ from $v_0$ to $v_0$, then go to the next vertex. Each time you enter a new vertex $v_i$, skip directly to the last visit of $v_i$ (bypassing another cycle), then keep going. (See Figure \ref{fig:to_rat_monoid}.)
\begin{center}
	\begin{tikzpicture}[scale=1.5]
	\draw[black, thick, -latex] (-.35,-.5)
		to [out=60, in=0, looseness=1.5] (-.4,-.3)
		to [out=180, in=140, looseness=1.5] (-.35,-.5)
		to [out=-40, in=-75, looseness=1.5] (0,-.1)
		to [out=105, in=-75, looseness=.9] (-.15,.75)
		to [out=105, in=160, looseness=1.5] (0,1.2)
		to [out=-20, in=30, looseness=1.5] (-.15,.75)
		to [out=-150, in=170, looseness=1.5] (0,-.1)
		to [out=-10, in=-35, looseness=2] (0.17,0.37)
		to [out=145, in=90, looseness=1.5] (-.25,0.4)
		to [out=-90, in=-120, looseness=2] (0.17,0.37)
		to [out=60, in=-120, looseness=1] (0.34,0.65);
	\node[black] at (.35,-.2) {$w$};
\end{tikzpicture}
	\hspace{10mm}
\begin{tikzpicture}[scale=1.5]
	\draw[black, thick, -latex] (-.35,-.5)
		to [out=-40, in=-75, looseness=1.5] (0,-.1)
		to [out=-10, in=-35, looseness=2] (0.17,0.37)
		to [out=60, in=-120, looseness=1] (0.34,0.65);
	\draw[LimeGreen, thick, -latex]	(-.35,-.5)
		to [out=60, in=0, looseness=1.5] (-.4,-.3)
		to [out=180, in=140, looseness=1.5] (-.35,-.5);
	\draw[Turquoise, thick, -latex]	(0,-.1)
		to [out=105, in=-75, looseness=.9] (-.15,.75)
		to [out=105, in=160, looseness=1.5] (0,1.2)
		to [out=-20, in=30, looseness=1.5] (-.15,.75)
		to [out=-150, in=170, looseness=1.5] (0,-.1);
	\draw[Magenta, thick, -latex] (0.17,0.37)
		to [out=145, in=90, looseness=1.5] (-.25,0.4)
		to [out=-90, in=-120, looseness=2] (0.17,0.37);
		
	\begin{footnotesize}
		\node[black] at (.7,-.2) {$s_1s_2\ldots s_\ell$};
		\node[LimeGreen] at (-.55,-.2) {$w_0$};
		\node[Turquoise] at (-.5,.7) {$w_1$};
		\node[magenta] at (.075,.6) {$w_2$};
	\end{footnotesize}	
\end{tikzpicture}

	\captionsetup{margin=6mm, font=footnotesize}
	
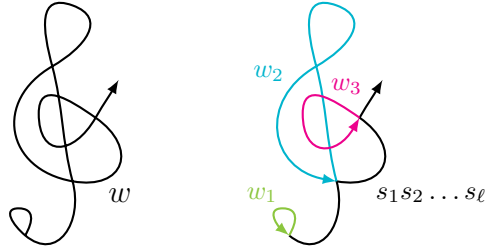
\captionof{figure}{A path in the automaton and its decomposition} \label{fig:to_rat_monoid}
\end{center}
$\blacktriangleright$ It follows that
\[ \Lc = \;\bigcup \; \Lc_{v_0\to v_0} \cdot s_1 \cdot \Lc_{v_1\to v_1} \cdot s_2 \cdot\ldots\cdot s_\ell \cdot \Lc_{v_\ell\to v_\ell} \]
where the union is taken over simple paths $v_0\overset{s_1}\to v_1\overset{s_2}\to \ldots \overset{s_\ell}\to v_\ell$ with $v_\ell\in\accept$, and $\Lc_{v\to v}$ is the language of words labeling cycles $v\to v$. The union is finite, and each language $\Lc_{v\to v}$ is regular (accepted by the automaton $M_{v\to v}=(V,\delta,v,\{v\})$).

\medbreak

If we managed to find \emph{bounded} regular languages $\Lc'_{v\to v}$ such that $\ev(\Lc'_{v\to v})=\ev(\Lc_{v\to v})$, we would be done with the language
\[ \Lc' = \bigcup_{s_1s_2\ldots s_\ell} \Lc'_{v_0\to v_0} \cdot s_1 \cdot \Lc'_{v_1\to v_1} \cdot s_2 \cdot\ldots\cdot s_\ell \cdot \Lc'_{v_\ell\to v_\ell} \]
which is bounded, regular, and evaluates to $R$. This is the subject of the next subsection.

\subsection{Main discussion} \label{ssec:cases}

We now fix $G=H_3(\Z)$, and consider a trim automaton $M=(V,\delta,v,\{v\})$ recognizing a language $\Lc$ and a rational subset $R$. In particular $\Lc$ is a submonoid of $G^\star$, and $R$ is a submonoid of $H_3(\Z)$. As in Proposition \ref{thm:Gil}, let
\[ X = \left\{ \bar t\, \bar u\, \bar t^{-1}\;\middle|\; \exists p\in V \text{ s.t.} \begin{array}{c} t\in G^\star \text{ labels a simple path }v\to p \\ u\in G^\star \text{ labels a simple cycle }p\to p \\ \text{both paths only intersect at }p \end{array} \right\}. \]
We will work with both $\Lc$ and $\tilde\Lc$ defined in Proposition \ref{thm:Gil}. (Note that $Y=\{\varepsilon\}$.)

\medbreak

Let $\pi\colon H_3(\Z)\onto\Z^2\colon g\mapsto \hat g$ be the abelianization map. We discuss depending on the subset positively spanned by $\pi(X)$, i.e., depending on $\{\lambda_1 y_1+\ldots+\lambda_r y_r\mid \lambda_i\in\R_{\ge 0},\, y_i\in \pi(X)\}$.
\begin{enumerate}[leftmargin=8mm, label=(\arabic*)]
	\item If $\pi(X)$ is included in a line.
	\item If $\pi(X)$ spans the whole plane (i.e., $\0$ belong to the interior of $\Conv(\pi(X))$).
	\item If $\pi(X)$ spans a half-plane.
	\item If $\pi(X)$ spans a cone. 
\end{enumerate}
In each case, we provide an (effectively computable) bounded regular language $\Lc'$ such that $\ev(\Lc')=\ev(\Lc)$. Case (4) will take most of our time.

\smallskip

\subsubsection{$\pi(X)$ spans $\{\0\}$, a ray or a line} \label{sssec:line}

The subgroup $\la X\ra$ is abelian (isomorphic to $\{e\}$, $\Z$ or $\Z^2$, but one might as well work in $\Z^X$). As $R$ is rational in $\la X\ra$ (Proposition \ref{thm:Gil}), we can compute a bounded regular language $\Lc'\subset \la X\ra^\star$ representing $R$ (Corollary \ref{cor:bounded_in_virtab}).

\smallskip

\subsubsection{$\pi(X)$ spans the whole plane}

For each $x=\overline t\overline u\overline t^{-1}\in X$ and $m\in\Z_{\ge 0}$, the set $R$ contains $x^m\ev(t w_{p\to v})=\ev(tu^mw_{p\to v})$, where $w_{p\to v}$ labels a path $p\to v$. Note that
\[
\pi(x^m\ev(tw_{p\to v})) = m\pi(x)+\pi(\ev(tw_{p\to v})) = m\pi(x) + O(1).
\]
It follows that $\Conv(\pi(R))=\R^2$, hence $R$ is a subgroup by Proposition \ref{prop:BCM}. More precisely $R=\la X\ra$ since $X\subseteq R\cdot R^{-1}$ (proof of Proposition \ref{thm:Gil}). We compute

\newpage

\begin{itemize}[leftmargin=8mm]
	\item $a,b\in \la X\ra$ such that $\pi(a),\pi(b)$ form a basis of $\pi(\la X\ra)\le \Z^2$. These can be found using Gaussian elimination on the matrix with vectors $\pi(x)$ ($x\in X$) as rows.
	\item $z^d\in \la X\ra$ such $z^d$ is a basis of $\la X\ra\cap[G,G]$. ($d$ is the smallest positive integer such that $z^d \in\la X\ra$. Check if $z,z^2,\ldots$ belong to $\la X\ra$, for instance using the solution to the Subgroup Membership \cite{Malcev}. The algorithm terminates since $z^{\det(\pi(a);\pi(b))}=[a,b]\in\la X\ra$.)
\end{itemize}
These two conditions mean that $\{a,b,z^d\}$ is a Mal'cev basis of $\la X\ra$. We deduce a regular normal form $\Lc'=\{a^pb^q(z^d)^r \mid p,q,r\in\Z\}\subset G^\star$ for $R$.

\medskip

\subsubsection{$\pi(X)$ spans a half-plane} \label{sssec:half}

By hypothesis, we can find
\begin{itemize}[leftmargin=7mm]
	\item $\bar s\,\bar a\,\bar s^{-1}$ and $\bar t\,\bar c\,\bar t^{-1}\in X$ such that $\hat a,\0$ and $\hat c\in\Z^2$ are on a line, in that order. We also fix $\tilde s\in G^\star$ labeling a path from the endpoint of $s$ to $v$. Define $\tilde t\in G^\star$ similarly.
	\item $b\in \Lc$ such that $\hat b$ doesn't lie on the same line. (Take $\bar r\,\bar u\bar r^{-1}\in X$ such that $\hat u$ doesn't lie on the line. Fix $\tilde r\in G^\star$ from the endpoint of $r$ to $v$, and take $b=ru\tilde r$.)
\end{itemize}
For each $\mathbf v\in\Z^2$, there exists unique $\alpha(\mathbf v)$, $\beta(\mathbf v)\in\Q$ such that $\mathbf v=\alpha(\mathbf v)\cdot \hat a+\beta(\mathbf v)\cdot \hat b$. 
\begin{lem} \label{lem:case_c}
	There exists a computable $K\ge 0$ such that, if $g\in R$ satisfies $\beta(\hat g)\ge K$, then $g$ can be written as
	\[ g = \ev\!\big( s a^{n_1} \tilde s \cdot b^m \cdot  t c^{n_2} \tilde t \cdot b^m \cdot  s a^{n_3} \tilde s \cdot b^{n_4}  \cdot w \big)  \]
	where $m$ is fixed, $n_1,n_2,n_3,n_4\in\Z_{\ge 0}$, and $w$ varies in a fixed finite subset of $\Lc$.
\end{lem}
\begin{proof} Let $h(n_1,n_2,n_3,n_4)= \ev\big(s a^{n_1} \tilde s \cdot b^m \cdot  t c^{n_2} \tilde t \cdot b^m \cdot  s a^{n_3} \tilde s \cdot b^{n_4}\big)$. We fix $p,q\in\Z_{>0}$ such that $p\hat a+q\hat c=\0$. In particular, we have
	\begin{equation}
		\hat h(n_1+p,n_2+q,n_3,n_4) = \hat h(n_1,n_2,n_3,n_4) = \hat h(n_1,n_2+q,n_3+p,n_4).
		\label{eq:proj33}
	\end{equation}
	Let us see how the area changes under the same transformations: \vspace*{2mm}
\begin{align*}
A\big(h(n_1+p,n_2+q,n_3,n_4)\big) - A\big(h(n_1,n_2,n_3,n_4)\big) & =  \underbrace{pA(a)+qA(c)+ \det\big(p\hat a;m\hat b+\pi(\tilde st)\big)}_{d_+(m)} \\
A\big(h(n_1,n_2+q,n_3+p,n_4)\big) - A\big(h(n_1,n_2,n_3,n_4)\big) & =  \underbrace{pA(a)+qA(c)- \det\big(p\hat a;m\hat b+\pi(\tilde ts)\big)}_{d_-(m)}
\end{align*}
Equation (\ref{eq:proj33}) and Lemma \ref{lem:compare_mod_N} imply that $d_+(m)$ and $d_-(m)$ are integers. We fix $m$ large enough so that $d_+(m)\cdot d_-(m)<0$ and fix $d=\gcd(d_+(m),d_-(m))$.
\medbreak
The subgroup $N\coloneqq\big\langle \bar a^d,\bar b^d, \bar c^d, z^d\big\rangle$ has finite index (Proposition \ref{prop:BCM}). Moreover
\[ \forall g\in N, g'\in G,\quad [g,g']=z^{\det(\hat g;\hat g')}\in\langle z^d\rangle \]
as $\hat g\in d\langle \hat a,\hat b,\hat c\rangle$. It follows that $[N,G]\le N$: the subgroup $N$ is normal. The quotient Cayley graph can be constructed using Lemma \ref{lem:compare_mod_N}. Using Lemma \ref{lem:FinInd_inter}, we can construct a finite subset $\Kc\subset\Lc$ of representatives for each coset of $N$ intersecting $R$.

We take $K=\beta(\pi(h(0,0,0,0)))+\max_{w\in\Kc}\beta(\hat w)$. Recall that $R$ is a submonoid in $G$, hence its image in the finite quotient $G/N$ is a subgroup. It follows that, for any $g\in R$, we can find $w\in\Kc$ such that $\bar w=  h(0,0,0,0)^{-1} g$ in $G/N$ (as $g, h(0,0,0,0)\in R$).

If furthermore $\beta(\hat g)\ge K$, there exists $m_1,m_2,m_4\in\Z_{\ge 0}$ such that
\begin{align*}
\pi(g) & = \pi\big(h(0,0,0,0)\big) + d\big(m_1\pi(\bar a)+m_2\pi(\bar c)+m_4\pi(\bar b)\big) + \pi(\bar w) \\
& = \pi\big(h(dm_1,dm_2,0,dm_4)\bar w\big)
\end{align*}
(We can ensure $m_1,m_2\ge 0$ as $\pi(\bar a),\pi(\bar c)$ are colinear, in opposite direction. Moreover, $dm_4=\beta(\hat g)-\beta(\pi(h(0,0,0,0)))-\beta(\hat w)\ge \beta(\hat g)-K\ge 0$.) Since $g$, $h(0,0,0,0)\bar w$ and $h(dm_1,dm_2,0,dm_3)\bar w$ lie in the same coset of $N$ and $\hat g=\hat h(dm_1,dm_2,0,dm_3)$, we have
\[ A(g) \equiv A\big(h(dm_1,dm_2,0,dm_4)\bar w\big) \pmod d \]
(Lemma \ref{lem:compare_mod_N}). Using the two transformations described above, we can find $n_1,n_2,n_3,n_4\in\Z_{\ge 0}$ such that $g=h(n_1,n_2,n_3,n_4)\bar w$.
\end{proof}

\medbreak

Finally, we can decompose into two (effectively computable) regular languages
\[ \Lc = \{w\in\Lc \mid \beta(\hat w)\ge K \} \sqcup \{ w\in \Lc \mid \beta(\hat w)<K\} =: \Lc_\reg\sqcup  \Lc_\abn. \]
\begin{itemize}[leftmargin=8mm]
	\item The first term can be replaced by a bounded regular language
	\[ \Lc'_\reg = \{ s a^{n_1} \tilde s \cdot b^m \cdot  t c^{n_2} \tilde t \cdot b^m \cdot  s a^{n_3} \tilde s \cdot b^{n_4}  \cdot w \mid w\in \Kc, n_1,n_2,n_3,n_4\in\Z_{\ge 0} \} \subseteq \Lc. \]
	Using Lemma \ref{lem:case_c}, we have $\ev(\Lc_\reg)\subseteq \ev(\Lc'_\reg)\subseteq \ev(\Lc)=R$.
	
	\item For the second term, we compute a trim automaton for $\Lc_\abn$ and compute the set $X$ associated to each strongly connected component. For each component, $\pi(X)$ is contained in the line through $\0$ and $\pi(\bar a)$. (Otherwise we would be able to pump along a cycle $u$ for which $\beta(\hat u)>0$ hence produce words $w\in\Lc_\abn$ with $\beta(\hat w)$ arbitrarily large.)
	
	We can therefore apply the arguments of \S\ref*{ssec:reduc} and \S\ref*{sssec:line} to get a bounded regular language $\Lc'_\abn$ such that $\ev(\Lc_\abn)=\ev(\Lc'_\abn)$. 
\end{itemize}
The language we are looking for is $\Lc'=\Lc'_\reg\cup \Lc'_\abn$.

\subsubsection{$\pi(X)$ spans a cone} \label{sssec:hard}
	
	We construct a bounded regular language $\Lc'_+$ such that
	\[ \{ g\in \ev(\Lc) \mid A(g)\ge 0\} \subseteq \ev(\Lc'_+) \subseteq \ev(\Lc). \]
	Obviously, we can do the same thing for elements of negative area, hence taking the union of both languages gives the desired bounded regular language $\Lc'$.
	
	\medbreak
	
Consider
\begin{itemize}[leftmargin=5mm]
	\item $tat^{-1}$ with $t$ a simple path from $v$ to some $p$, and $a$ a simple loop from $p$ to $p$ such that $\pi(\bar a)$ belongs to the lower side of the cone, and $\frac{A(tat^{-1})}{\norm{\pi( a)}}$ is maximized, where $\norm{\vardot}$ is your favorite norm on $\R^2$. We also fix $\tilde t\in G^\star$ labeling a path back from $p$ to $v$.
		
	\item $s^{-1}bs$ with $s$ a simple path from some $q$ to $v$, and $b$ a simple loop from $q$ to $q$ such that $\pi(\bar b)$ belongs to the upper side of the cone, and $\frac{A(s^{-1}bs)}{\norm{\pi( b)}}$ is maximized. We also fix $\tilde s\in G^\star$ labeling a path from $v$ to $q$.
\end{itemize}
For each $\mathbf v\in\Z^2$, there exists unique $\alpha(\mathbf v)$, $\beta(\mathbf v)\in\Q$ such that $\mathbf v=\alpha(\mathbf v)\cdot \hat a+\beta(\mathbf v)\cdot \hat b$. Moreover, $\mathbf v$ belongs to the cone if and only if $\alpha(\mathbf v)\ge 0$ and $\beta(\mathbf v)\ge 0$.

\bigskip

\thesubsubsection.i. \emph{There does not exist $x\in X$ such that $\pi(x)=\mathbf 0$ and $A(x)>0$.}

This is the difficult case. If $\hat g$ is far from the border of the cone, and $A(g)$ is far from its maximum (for a fixed value of $\hat g$), then there is enough leeway to find a word representing $g$ in a bounded sub-language of $\Lc$ which we are going to specify.
	\begin{lem}\label{lem:conei}
		There exist computable $K,m\ge 0$ such that, if $g\in R$ satisfies
		\[ \begin{array}{c}
		0 \le A(g) \le \frac12\det(\hat a;\hat b) \cdot \alpha(\hat g)\beta(\hat g) + A(tat^{-1}) \cdot \alpha(\hat g) + A(s^{-1}bs) \cdot \beta(\hat g) - K, \\[2mm]
		\alpha(\hat g) \ge K \quad\text{\normalfont and}\quad \beta(\hat g)\ge K,
		\end{array} \]
		then $g$ can be written as
		\[ g= \ev\!\big(t a^{n_1} \tilde t \cdot \tilde s b^{n_2}s \cdot t a^{n_3}\tilde t \cdot w \cdot \tilde s b^{p_1}s \cdot t a\tilde t \cdot \tilde s b^{p_2}s \cdot \ldots \cdot t a\tilde t \cdot \tilde s b^{p_m}s\big) \]
		where $n_1,n_2,n_3,p_1,p_2,\ldots,p_m$ vary in $\Z_{\ge 0}$ and $w$ varies in a fixed finite subset of $\Lc$. 
	\end{lem}
\begin{proof}Let 
	\[ h(n_1,n_2,n_3;w;p_1,\ldots,p_m) = \ev\!\big(t a^{n_1} \tilde t \cdot \tilde s b^{n_2}s \cdot t a^{n_3}\tilde t \cdot w \cdot \tilde s b^{p_1}s \cdot t a\tilde t \cdot \tilde s b^{p_2}s \cdot \ldots \cdot t a\tilde t \cdot \tilde s b^{p_m}s\big). \]
	The operation $(p_i,p_{i+1})\to (p_i+1,p_{i+1}-1)$ preserves $\hat h$, and decreases the area by
	\[  A\big(h(\ldots,p_i,p_{i+1},\ldots)\big) - A\big(h(\ldots,p_i+1,p_{i+1}-1,\ldots)\big) = \det\!\big(\pi(s\cdot  ta\tilde t \cdot \tilde s);\hat b\big).\]
	Note that $st\tilde t\tilde s$ labels a cycle in the automaton, so $\alpha(\pi(st\tilde t\tilde s))\ge 0$. Let 
	\[ d= \det\big(\pi(s \cdot ta\tilde t \cdot \tilde s);\hat b\big)= \big(1+\alpha(\pi(st\tilde t\tilde s))\big) \cdot \det\!\big(\hat a;\hat b\big)\ge \det\!\big( \hat a;\hat b\big)>0. \]
	We also fix $m=2d+1$. We consider $N=\langle \bar a^d,\bar b^d,z^d\rangle$. Using Lemma \ref{lem:FinInd_inter}(c), we compute a finite set $\Kc\subset\Lc$ containing a representative for each coset of $N$ intersecting $R$.
	
	\medbreak
	
	For each $g\in R$ satisfying the condition from Lemma \ref{lem:conei}, we find parameters so that $h(n_1,n_2,n_3;w;p_1,\ldots,p_m)=g$. We split this into four steps:
	\begin{enumerate}[leftmargin=10mm, label=(\arabic*)]
		\item We find $w\in\Kc$ such that $g=h(0,0,0;w;0,\ldots,0)$ in $G/N$.
		\item For $K$ computably large enough, we find $N_1,P_m\ge 0$ such that
		\[ g=h(N_1,0,0;w;0,\ldots,0,P_m) \]
		in $G/(N\cap [G,G])$. Moreover $A\big(h(N_1,0,0;w;0,\ldots,0,P_m)\big)-A(g)\ge 0$.
		\item For $K$ large enough, we find $n_1,n_2,n_3,p_m\ge 0$ with $n_2\le \frac35\beta(\hat g)$ and
		\[ g=h(n_1,n_2,n_3;w;0,\ldots,0,p_m) \]
		in $G/(N\cap[G,G])$. Moreover
		\[ 0 \le  A(h(n_1,n_2,n_3;w;0,\ldots,0,p_m)) - A(g) \le d\left(\frac35 \beta(\hat g)+\beta(\pi(\tilde t\tilde sst))\right) \det(\hat a;\hat b). \]
		\item For $K$ large enough, we find $p_1,\ldots,p_m\ge 0$ such that $g=h(n_1,n_2,n_3;w;p_1,\ldots,p_m)$.
	\end{enumerate}

	\medbreak
	
	
	(1) Since the image of the monoid $R$ in $G/N$ is a subgroup, we can find $w\in\Kc$ such that
	\[ w = (t\tilde t \cdot \tilde ss \cdot t\tilde t)^{-1} \cdot g \cdot \big( (\tilde s s\cdot ta\tilde t)^{m-1} \cdot \tilde s s\big)^{-1} \qquad\text{in }G/N. \]
	It follows that $g=h(0,0,0;w;0,\ldots,0)$ in $G/N$.
	
	\medbreak
	
	(2) In particular, we have $\hat g-\hat h(0,0,0;w;0,\ldots,0)\in d\langle \hat a,\hat b\rangle$ (Lemma \ref{lem:compare_mod_N}). It follows that
	\begin{align*}
		N_1 & =\alpha(\hat g)-\alpha(\hat h(0,0,0;w;0,\ldots,0)), \\
		P_m & = \beta(\hat g)-\beta(\hat h(0,0,0;w;0,\ldots,0))
	\end{align*}
	are integer multiples of $d$. Moreover $N_1,P_m\ge 0$ for $K$ computably large enough. We have
	\[
	g= h(N_1,0,0;w;0,\ldots,0,P_m) \qquad\text{in }G/(N\cap[G,G]).
	\]
	The area of $h(N_1,0,0;w;0,\ldots,0,P_m)$ is given by
	\[ A(ta^{N_1} \cdot v \cdot b^{P_m}s) = \\ \frac12\det(\hat a;\hat b)\cdot \alpha(\hat g)\beta(\hat g) + A(tat^{-1}) \cdot \alpha(\hat g) + A(s^{-1}bs) \cdot \beta(\hat g) - C(w) \]
	with $v=\tilde t\tilde sst\tilde tw(\tilde ss ta\tilde t)^{m-1}\tilde s$, and $C(w)$ is a computable constant that only depends on $w$. For $K$ large enough, we may therefore assume $A(g)\le A(h(N_1,0,0;w;0,\ldots,0,P_m))$.
	\begin{center}
		\begin{tikzpicture}[scale=.5, rotate=-25]
	\newcommand{\asqui}[2]{
		\begin{scope}[shift={(#1,#2)}]
			\draw[-latex, looseness=4, out=30, in=210] (0,0) to (3,1);
			\node at (1.5,-.3) {$a$};
	\end{scope}}
	
	\newcommand{\tsqui}[2]{
		\begin{scope}[shift={(#1,#2)}]
			\draw[dotted, looseness=2, out=-20, in=170] (-1,1) to (0,0);
	\end{scope}}
	
	\newcommand{\bsqui}[2]{
		\begin{scope}[shift={(#1,#2)}]
			\draw[-latex, looseness=1.5, out=75, in=285] (0,0) to (0,2);
			\node at (0.6,.9) {$b$};
	\end{scope}}
	
	\newcommand{\ssqui}[2]{
		\begin{scope}[shift={(#1,#2)}]
			\draw[dotted, looseness=3, out=150, in=-30] (0,0) to (-1,1);
	\end{scope}}
	
	
	\fill[LimeGreen!25]
	(-1,1) to[looseness=2, out=-20, in=170]
	(0,0) to[looseness=4, out=30, in=210]
	(3,1) to[looseness=4, out=30, in=210]
	(6,2) to[looseness=4, out=30, in=210]
	(9,3) to[looseness=1.5, out=95, in=-180]
	(11,5) to[looseness=1.5, out=75, in=285]
	(11,7) to[looseness=1.5, out=75, in=285]
	(11,9) to[looseness=1.5, out=75, in=285]
	(11,11) to[looseness=1.5, out=75, in=285]
	(11,13) to[-latex, looseness=3, out=150, in=-30]
	(10,14) to (-1,1);
	
	\draw[dashed] (-1,1) -- (10,{14/3}) -- (10,14);
	
	\draw[-latex, looseness=2, out=-20, in=170] (-1,1) to (0,0);
	\asqui00
	\asqui31
	\asqui62
	
	\tsqui31
	\tsqui62
	\tsqui93
	
	\draw[-latex, looseness=1.5, out=95, in=-180] (9,3) to (11,5);
	
	\ssqui{11}5
	\ssqui{11}7
	\ssqui{11}9
	\ssqui{11}{11}
	\bsqui{11}5
	\bsqui{11}7
	\bsqui{11}9
	\bsqui{11}{11}
	\draw[-latex, looseness=3, out=150, in=-30] (11,13) to (10,14);

	\draw[fill=black] (-1,1) circle (3pt);
	\footnotesize{\node[circle, fill=black, inner sep=1.5pt, label=above:$\hat g$] at (10,14) {};}
	\node at (-1,.4) {$t$};
	\node at (9.8,4) {$v$};
	\node at (11,13.8) {$s$};
\end{tikzpicture}
		
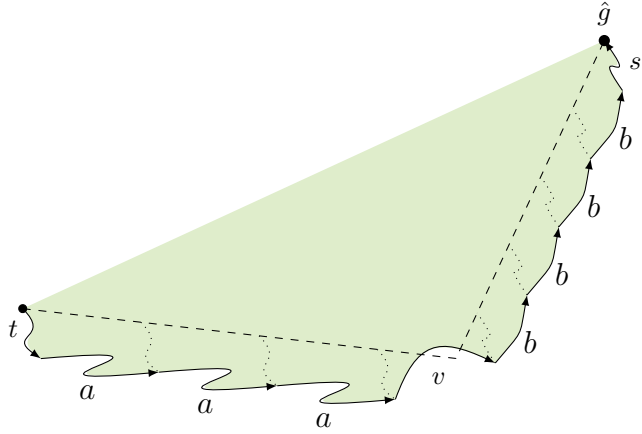
\captionof{figure}{The path corresponding to $h(N_1,0,0;w;0,\ldots,0,P_m)$.} \label{fig:triangle}
	\end{center}
	Indeed, the large triangle (Figure \ref{fig:triangle}) has area $\frac12\det(\hat a;\hat b)\cdot \alpha(\hat g)\beta(\hat g)$. The $N_1$ (resp.\ $P_m$) smaller regions bordered by $a$ (resp.\ $b$) have area $A(tat^{-1})$ (resp.\ $A(s^{-1}bs)$).
	
	\medbreak
	
	(3) Next, we apply two transformations in order:
	\begin{itemize}[leftmargin=8mm]
		\item $(n_2,p_m)\to (n_2+d,p_m-d)$. This doesn't change $\hat h$, and decreases the area by
		\begin{align*}
			& A\big(h(N_1,n_2,0;w;0,\ldots,0,p_m)\big) - A\big(h(N_1,n_2+d,0;w;0,\ldots,0,p_m-d)\big)  \\
			& = \det\!\Big(\pi\big(s \cdot t\tilde t \cdot w \cdot (\tilde ss \cdot ta\tilde t)^{m-1}\cdot \tilde s\big);d\cdot \pi(b)\Big) = d\Big((m-1)+\alpha\big(\pi(w(\tilde sst\tilde t)^m)\big)\Big)\det(\hat a;\hat b)> 0.
		\end{align*}
		(Note that $\alpha(\pi(\tilde ss)),\alpha(\pi(t\tilde t)),\alpha(\hat w)\ge 0$ since $\tilde ss$, $t\tilde t$, $w\in \Lc$.) Repeat as long as possible while keeping $n_2\le\frac35\beta(\hat g)$ and $A(h)\ge A(g)$. (The first condition ensures that $p_m\ge 0$ since $p_m=P_m-n_2\ge \frac25\beta(\hat g)-\beta(\hat h(0,0,0;w;0,\ldots,0))>0$ for $K$ large enough.)
		\item $(n_1,n_3)\to (n_1-d,n_3+d)$. This doesn't change $\hat h$, and decreases the area by
		\begin{align*}
		& A\big(h(n_1,n_2,n_3;w;0,\ldots,0,p_m)\big) - A\big(h(n_1-d,n_2,n_3+d;w;0,\ldots,0,p_m)\big) \\
		& = \det\!\big(d\pi(a);\pi(\tilde t \cdot \tilde sb^{n_2}s\cdot t)\big) = d\big(n_2+\beta(\pi(\tilde t\tilde ss t))\big)\det(\hat a;\hat b)\ge 0.
		\end{align*}
		Repeat as long as possible while keeping $n_1\ge 0$ and $A(h)\ge A(g)$.
	\end{itemize}
	At the end of this process, we have found an element $h(n_1,n_2,n_3;w;0,\ldots,0,p_m)$ with $n_1,n_2,n_3,p_m\ge 0$ and $n_2\le\frac 35\beta(\hat g)$, which coincides with $g$ in $G/(N\cap[G,G])$. In particular, $A(h)-A(g)\in d\Z$. Moreover, since we cannot apply the second step anymore, we have
	\[ 0\le A(h)-A(g) \le d\big(n_2+\beta(\pi(\tilde t\tilde ss t))\big)\det(\hat a;\hat b) \le d\left(\frac35 \beta(\hat g)+\beta(\pi(\tilde t\tilde sst))\right) \det(\hat a;\hat b). \]
	Here we use that $A(g)\ge 0$. If we had continued until the conditions $n_2\le\frac35\beta(\hat g)$ and $n_1\ge 0$ were the limiting conditions, we would have reached an area of
	\[ A(h) \sim A\left(h\left(0,\frac35P_m,N_1;w;0,\ldots,0,\frac25P_m\right)\right) \sim -\frac1{10}\det(\hat a;\hat b)\cdot \alpha(\hat g)\beta(\hat g)<0\le A(g) \]
	where $\sim$ means that the quotient converges to $1$ when $\min\{\alpha(\hat g),\beta(\hat g)\}\to\infty$. For $K$ computably large enough, this implies that the condition $A(h)\ge A(g)$ stopped us.
	
	\medbreak
	
	(4) Finally, we use the operation $(p_i,p_{i+1})\to (p_i+1,p_{i+1}-1)$. We can use this operation up to $(m-1)(P_m-n_2) = 2d(P_m-n_2)$ times, decreasing the area by exactly $d$ each time. We can therefore reduce the area by any multiple of $d$ up to
	\[ 2d\cdot \left(\frac25\beta(\hat g)-\beta(\hat h(0,0,0;w;0,\ldots,0))\right)\cdot d \ge d\left(\frac35 \beta(\hat g)+\beta(\pi(\tilde t\tilde sst))\right) \det(\hat a;\hat b) \]
	as $d\ge \det(\hat a;\hat b)$ and $\beta(\hat g)\ge K$ with $K$ large enough. Since $A(h)-A(g)$ was an integer multiple of $d$ in the correct range after step (3), we get $A(h)=A(g)$ at some point.
\end{proof}
	
\bigskip

As in the Case \ref{sssec:half}, we decompose $\Lc$ into two regular languages, with all elements in the first set treated by Lemma \ref{lem:conei}, and all elements in the second case following the border of the cone (hence being easier to treat). We first need a lemma to take care of the \say{area condition}  with a finite state automaton.	
\begin{lem}
	Consider a word $w\in\Lc$ and $\tilde w=x_1x_2\ldots x_\ell\in X^*$ its decomposition. Let
	\begin{align*}
	M(w)\, & = \# \left\{ i \;\big|\; 1\le i\le \ell \text{ and } \alpha(x_i),\beta(x_i)>0 \right\}, \\
	N(w)\, & = \# \left\{ (i,j) \;\big|\; 1\le i<j\le \ell \text{ and } \alpha(x_j),\beta(x_i)>0 \right\}.
	\end{align*}
	There exists a computable constant $L$ such that, if
	\[ A(\bar w) \ge \frac12\det(\hat a;\hat b)\cdot \alpha(\hat w)\beta(\hat w) + A(tat^{-1}) \cdot \alpha(\hat w) + A(s^{-1}bs) \cdot \beta(\hat w) - K, \]
	then $M(w),N(w)\le L$.
	\end{lem}
\begin{proof}
	Let
	\begin{align*}
		\Delta & = \max\Big(\big\{ A(x) \mid x\in X \big\}\cup\{0\}\Big), \\
		\delta & = \min\Big(\big\{ \alpha(x),\beta(x)\mid x\in X\big\} \setminus\{0\}\Big).
	\end{align*}
	We can bound the area of $\bar w$ by
	\begin{align*}
	A(\bar w) & = \sum_{i<j} \frac12 \det(\hat x_i;\hat x_j) + \sum_i A(x_i) \\
	& = \frac12\det(\hat a;\hat b) \left(\sum_{i<j} \alpha(\hat x_i)\beta(\hat x_j)-\alpha(\hat x_j)\beta(\hat x_i)\right) + \sum_i A(x_i) \\
	& \le \frac12\det(\hat a;\hat b) \left(\sum_{i} \alpha(\hat x_i) \cdot \sum_j\beta(\hat x_j)-N(w)\cdot \delta^2\right) + \\
	& \quad +\sum_{i:\beta(\hat x_i)=0} \alpha(\hat x_i)A(tat^{-1}) + \sum_{i:\alpha(\hat x_i)=0} \beta(\hat x_i)A(s^{-1}bs) + M(w)\cdot \Delta \\
	& \le \frac12\det(\hat a;\hat b)\cdot \alpha(\hat w)\beta(\hat w) + A(tat^{-1}) \cdot \alpha(\hat w) + A(s^{-1}bs) \cdot \beta(\hat w)\; + \\[1mm] & \quad + M(w)\cdot \Delta - \frac12\det(\hat a;\hat b)\cdot N(w)\cdot \delta^2.
	\end{align*}
	The first inequality is the only part using that $A(x)\le 0$ for all $x\in X$ such that $\pi(x)=\0$, and that $A(tat^{-1})/\norm{\pi(a)}$ and $A(s^{-1}bs)/\norm{\pi(b)}$ are maximal. This ensures that
	\begin{itemize}[leftmargin=8mm]
		\item $A(x_i)\le 0$ if $\alpha(\hat x_i)=\beta(\hat x_i)=0$,
		\item $A(x_i)\le \alpha(x_i)A(tat^{-1})$ if $\alpha(\hat x_i)\ne 0$ and $\beta(\hat x_i)=0$,
		\item $A(x_i)\le \beta(x_i)A(s^{-1}bs)$ if $\alpha(\hat x_i)= 0$ and $\beta(\hat x_i)\ne0$.
	\end{itemize} 
	It follows that $M(w)\cdot \Delta - \frac12\det(\hat a;\hat b) \delta^2 \cdot N(w) \ge -K$. Since ${M(w)\choose 2}\le N(w)$ we get
	\[ K \ge \frac{\delta^2\det(\hat a;\hat b)}4 \cdot M(w)\big(M(w)-1\big) - \Delta \cdot M(w) \]
	hence $M(w)\le L$ for some computable $L$. Finally $N(w)\le \frac2{\delta^2\det(\hat a;\hat b)}(M(w)\cdot \Delta+K)$.
\end{proof}

\medbreak

We decompose the language $\tilde\Lc$ given by Proposition \ref{thm:Gil} into two regular subsets:
\[
\Lc_\reg\coloneqq \left\{ \tilde w\in\tilde\Lc \;\middle|\; \alpha(\hat w),\beta(\hat w)\ge K \text{ and } N(w)\ge L \right\},
\]
and $\Lc_\abn\coloneqq \tilde\Lc\setminus\Lc_\reg$.
\begin{itemize}[leftmargin=8mm]
	\item We replace $\Lc_\reg$ by the bounded regular language
	\[ \Lc'_\reg = \bigcup_{w\in\Kc} t a^* \tilde t \cdot \tilde s b^*s \cdot t a^*\tilde t \cdot w \cdot (\tilde s b^*s \cdot t a\tilde t)^{m-1} \cdot \tilde s b^*s, \]
	so that $\{ g\in \ev(\Lc_\reg) \mid A(g)\ge 0\}  \subseteq \ev(\Lc'_\reg)\subseteq R$.
	
	\item For the second term, we compute a trim automaton for $\Lc_\abn$ and compute the set $X$ associated to each strongly connected component. For each component, $\pi(X)$ is contained in the line through $\pi(a)$ or through $\pi(b)$. We can apply \S\ref*{ssec:reduc} and \S\ref*{sssec:line} to get a bounded regular language $\Lc'_\abn$ such that $\ev(\Lc'_\abn)=\ev(\Lc_\abn)$. 
\end{itemize}
The language we are looking for is $\Lc_+'=\Lc'_\reg\cup \Lc'_\abn$.

\newpage

\thesubsubsection.ii. \emph{There does exist $x\in X$ such that $\pi(x)=\mathbf 0$ and $A(x)>0$.}

This case is similar to \S\ref*{sssec:half}. We fix $r,c,\tilde r\in G^\star$ labeling a path $v\to p$, a cycle $p\to p$ and a path $p\to v$ respectively, such that $\hat c=\mathbf 0$ and $A(\bar c)>0$.
\begin{lem} \label{lem:coneii}
	There exists a computable $K\ge 0$ such that, if $g\in R$ satisfies
	\[ A(g) \ge 0, \quad \alpha(\hat g) \ge K \quad\text{\normalfont and}\quad \beta(\hat g)\ge K, \]
	then $g$ can be written as $g=\ev\!\big(\tilde sb^{n_1}s\cdot ta^{n_2}\tilde t\cdot rc^{n_3}\tilde r \cdot w \big)$ where $n_1,n_2,n_3$ varies in $\Z_{\ge 0}$ and $w$ varies in a fixed finite subset of $\Lc$.
\end{lem}
\begin{proof} Let $d=A(c)>0$. We introduce the notation
	\[ h(n_1,n_2,n_3) = \ev\!\big(\tilde sb^{n_1}s\cdot ta^{n_2}\tilde t\cdot rc^{n_3}\tilde r \big).\]
	We consider the finite-index normal subgroup $N=\langle \bar a^d,\bar b^d,\bar c=z^d\rangle$. Using Lemma \ref{lem:FinInd_inter}(c), we construct a finite subset $\Kc\subset \Lc$ of representatives for each coset of $N$ intersecting $R$. For $g\in R$ satisfying the condition from the Lemma \ref{lem:coneii}, we take
	\begin{itemize}[leftmargin=8mm]
		\item $w\in\Kc$ such that $\bar w = h(0,0,0)^{-1}g$ in $G/N$.
		\item $n_1=\alpha(\pi(g))-\alpha(\pi(h(0,0,0)\bar w))$ and $n_2=\beta(\hat g)-\beta(\pi(h(0,0,0)\bar w))$. These are integer multiples of $d$ (Lemma \ref{lem:compare_mod_N}). Moreover, they are non-negative for $K$ large enough.
		\item $n_3=\frac1d\left(A(g)-A\big(h(n_1,n_2,0)\cdot \bar w\big)\right)$. Note that $n_3$ is an integer since
		\[ g=h(n_1,n_2,0)\bar w \qquad\text{in }G/(N\cap[G,G]). \]
		Moreover, $n_3\ge 0$ for $K$ large enough. Indeed, we have $n_1,n_2\ge K-C_1(w)$ hence
		\[ A\big(h(n_1,n_2,0)\cdot \bar w\big) = -\frac12\det(\hat a;\hat b)\cdot n_1n_2 + C_2(w)\cdot n_1 + C_3(w)\cdot n_2 + A(w) \le 0 \le A(g)\]
		for $K$ large enough, where $C_1(w),C_2(w),C_3(w)$ are constant depending on $w$.
	\end{itemize}
	 We conclude with $\pi(h(n_1,n_2,n_3)\bar w)=\pi(h(n_1,n_2,0)\bar w)=\pi(g)$ and $A(h(n_1,n_2,n_3)\bar w)=A(h(n_1,n_2,0)\bar w)+n_3d=A(g)$, hence $h(n_1,n_2,n_3)\bar w=g$.
\end{proof}

\bigskip

Finally, we decompose $\Lc$ into two regular languages
\[ \Lc = \{w\in \Lc \mid \alpha(\hat w),\beta(\hat w)\ge K\} \sqcup \{w\in \Lc \mid \alpha(\hat w)< K \text{ or }\beta(\hat w)< K\} =: \Lc_{\reg}\sqcup \Lc_{\abn}. \]
\begin{itemize}[leftmargin=8mm]
	\item The first term can be replaced by a bounded regular language
	\[ \Lc'_{\reg} = \left\{ \tilde s b^{n_1} s \cdot ta^{n_2}\tilde t \cdot rc^{n_3}\tilde r \cdot w \;\big|\; n_1,n_2,n_3\ge 0, w\in\Kc\right\} \subseteq \Lc. \]
	Using Lemma \ref{lem:coneii}, we have $\left\{g\in \ev(\Lc_\reg) \mid A(g)\ge 0\right\} \subseteq \ev(\Lc'_\reg) \subseteq \ev(\Lc)$.
	\item For the second term, we compute a trim automaton for $\Lc_\abn$ and compute the set $X$ associated to each strongly connected component. For each component, $\pi(X)$ is contained in the line through $\pi(a)$ or through $\pi(b)$. We can apply \S\ref*{ssec:reduc} and \S\ref*{sssec:line} to get a bounded regular language $\Lc'_\abn$ such that $\ev(\Lc'_\abn)=\ev(\Lc_\abn)$. 
\end{itemize}
The language we are looking for is $\Lc_+'=\Lc'_\reg\cup \Lc'_\abn$.

\subsection{Equations under rational constraints} \label{ssec:thmC}

Finally, we prove that $\RSMP(H_3(\Z))$ is decidable, using Proposition \ref{thm:reduc}. We prove a stronger result about solving equations under rational constraints.
\begin{defi}
	An equation with rational constraints in $G$ consists of
	\begin{itemize}[leftmargin=8mm]
		\item an element $w\in F(x_1,\ldots,x_n)*G$ and
		\item $n$ rational subsets $R_1,\ldots,R_n\subseteq G$.
	\end{itemize}
	This equation admits a solution if there exists $(g_1,\ldots,g_n)\in R_1\times \ldots \times R_n$ such that $w(g_1,\ldots,g_n):=f(w)=e_G$, where $f$ is the homomorphism
	\[ f\colon \begin{pmatrix}
		F(x_1,\ldots,x_n) * G & \longto & G \\ g & \longmapsto & g \\ x_i & \longmapsto & g_i
	\end{pmatrix}.\]
\end{defi}
This generalizes the Rational Subset Membership $g\overset?\in R$, as we may consider the equation $x_1=g$ under the rational constraint $x_1\in R$.
\begin{thm} \label{thm:rat_memb}
	There exists an algorithm which takes as input an equation with rational constraints in $H_3(\Z)$, and decides whether it admits a solution.
\end{thm}
This extends the analogous result without rational constraints, due to Duchin, Liang and Shapiro. The proof is a straightforward adaptation of \cite[Theorem 3]{duchin2015equations} and \cite[Theorem 6.8]{Knapsack_Heis}, with Proposition \ref{thm:reduc} as a starting point.
\begin{proof}
	We first prove the statement under the extra assumption that each rational constraint is given as $R_i = h_{i,0}\, \{k_{i,1}\}^*\, h_{i,1} \, \{k_{i,2}\}^*\, h_{i,2} \ldots h_{i,\ell_i-1}\, \{k_{i,\ell_i}\}^* \, h_{i,\ell_i}$. (Such a set will be called \say{Knapsack-like}.) In particular, a generic element of $R_i$ is given as
	\[ g_i=g_i(n_1,\ldots,n_{\ell_i})=h_{i,0}\,k_{i,1}^{n_1}\,h_{i,1}\,k_{i,2}^{n_2}\,h_{i,2}\ldots h_{i,\ell_i-1}\,k_{i,\ell_i}^{n_{\ell_i}}\,h_{i,\ell_i}\]
	with $n_j\in\N$. In coordinates, the existence of a solution reduces to the system
	\[ \begin{cases}
		\hspace*{4.5mm} \hat w(g_1,\ldots,g_n) \hspace*{1.4mm} = \mathbf 0 \\ A(w(g_1,\ldots,g_n))  = 0 
	\end{cases}\]
	This system consists of two linear and one quadratic equations, with coefficients in $\frac12\Z$ and unknowns $n_{i,j}\in\N$, hence we can decide if it admits a solution using \cite{Grunewald_Segal}.
	
	In general, we can write each rational constraint as $R_i = \bigcup_{j=1}^{m_i} R_{ij}$ where each $R_{ij}$ is Knapsack-like using Proposition \ref{thm:reduc}. Therefore, we only have to check whether one of $m_1\ldots m_n$ systems with Knapsack-like constraints admits a solution. The final answer is \say{Yes} if any of these answers is yes, \say{No} otherwise.
\end{proof}

\newpage

\begin{rem}
	It should be noted that the result about Rational Subset Membership can be extended to $H_3(\Q)$ as every finitely generated subgroup of $H_3(\Q)$ is (effectively) isomorphic to a subgroup of $H_3(\Z)$. Indeed, given a finite set of matrices
	\[ M_i = \begin{pmatrix} 1 & a_i & c_i \\ 0 & 1 & b_i \\ 0 & 0 & 1 \end{pmatrix} \in H_3(\Q)\quad (i=1,2,\ldots,r), \]
	find $N\ne 0$ such that $Na_i, Nb_i, N^2c_i\in\Z$. The dilation $\delta_N\colon (a,b,c)\mapsto (Na,Nb,N^2c)$ is an automorphism $H_3(\Q)\to H_3(\Q)$ such that $\delta_N\big(\!\la M_i\ra\!\big)\le H_3(\Z)$.
\end{rem}
\section{Further questions and remarks}  \label{sec:Qs}

One may wonder for which class of groups can the Rational Subset Membership be reduced to the Knapsack problem (as in Proposition \ref{thm:reduc}):
\begin{adjustwidth}{2.5mm}{2.5mm}
\textbf{Question D.} Characterize groups $G$ such that every rational subset $R\subseteq G$ can be effectively represented by a \emph{bounded} regular language $\Lc'\subset G^\star$. In particular, what if
\begin{enumerate}[leftmargin=8mm]
	\item $G$ is $2$-step nilpotent with cyclic derived subgroup?
	\item $G$ is $2$-step nilpotent?
	\item $G=E$ is the following \say{filiform} $3$-step nilpotent group
	\[ E = \la x,y_1,y_2 \;\big|\; [x,y_i]=y_{i+1}\text{ for }i=1,2,\; [x,y_3] =[y_i,y_j]=1\ra\;\text? \vspace*{-2mm} \]
	\item Does any group of super-polynomial growth have this property? \vspace*{2mm}
\end{enumerate}
\end{adjustwidth}
$\blacktriangleright$ Using Proposition \ref{thm:Gil} and Lemma \ref{lem:FinInd_inter}, we can prove that this property passes to quotients, subgroups and finite-index overgroups.

\medbreak

$\blacktriangleright$ Recall that $2$-step nilpotent groups with infinite cyclic derived subgroup contain copies of $H_{2m+1}(\Z)\times\Z^n$ as finite-index subgroups, see \cite[Lemma 7.1]{Stoll_series}. In particular, we only have to consider $G=H_{2m+1}(\Z)\times\Z^n$ to answer Question D.1.

A positive answer to Question D.1 would be particularly interesting as it would imply decidability for the Rational Subset Membership for those groups. Indeed, the Knapsack problem is also decidable for those groups (adapting the proof of \cite[Theorem 6.8]{Knapsack_Heis} from $H_3(\Z)\times\Z^n$ to $H_{2m+1}(\Z)\times\Z^n$, then using \cite[Theorem 7.3]{Knapsack_Heis}).


\medbreak

$\blacktriangleright$ An instance of Question D.2 of special interest is
\[ G = N_{2,r} = \la x_1,x_2,\ldots,x_r \;\big|\; \text{2-step nilpotent} \ra \]
with $R=\{x_1,x_2,\ldots,x_r\}^*$. A positive answer (algorithmic in $r$) would allow to reduce the Submonoid Membership Problem to the Knapsack Problem in any $2$-step nilpotent groups. In turn, this would provide an example of a group with decidable Submonoid Membership and undecidable Knapsack problem (adapting the proof of Corollary \ref{cor:exe}). Note that the answer is positive for $r=2$, as we have the equality
\[ \{x,y\}^* = \{x\}^*\{y\}^*\,x\,\{y\}^*\{x\}^* \,\sqcup\, \{y\}^* \]
in $N_{2,2}=H_3(\Z)$. This equality is one of the key motivations behind \S\ref{sssec:hard}.
\medbreak

$\blacktriangleright$ The property fails for the free $3$-step nilpotent group of rank $2$
\[ G = N_{3,2} = \la x,y \;\big|\; \text{3-step nilpotent} \ra. \]
This group admits a path model similar to the model for elements of $H_3(\Z)$ defined in \S3.1, where we keep track of an additional parameter $B(g)\in\R^2$ (see \cite[Section 3.1]{horo} for all necessary definitions). We argue that $R=\{x,y\}^*$ cannot be represented by a bounded regular language. Let $\Lc'\subset G^\star$ be a bounded regular language such that $\ev(\Lc')\subseteq R$, say
\[ \Lc' = \bigcup_{i=1}^I v_{i,0}\, \{w_{i,1}\}^*\, v_{i,1}\, \{w_{i,2}\}^*\, v_{i,2} \ldots v_{i,\ell_i-1}\{w_{i,\ell_i}\}^*\, v_{i,\ell_i}.\]
Observe that $\hat w_{i,j}\in\Z_{\ge 0}^2$, otherwise the word $w=v_{i,0}\,v_{i,1}\ldots v_{i,j-1}\,w_{i,j}^N\,v_{i,j}\ldots v_{i,\ell_i}\in\Lc'$ will project to $\hat w=N\cdot\hat w_{i,j}+O(1)\notin\Z_{\ge 0}^2=\pi(R)$ for $N$ large enough.

Take a direction $u\in\R_{\ge 0}^2$ which is not proportional to any of the $\hat w_{i,j}$ and take an element $\gamma_{u,n}\in\ev\{x,y\}^*$ of length $n$ which best follows the ray $\R^+u$, as defined in \cite[Section 3.2]{horo}. Adapting computations from \cite[Section 4]{Intermediate}, we have
\[ \la B(\gamma_{u,n})-B(\bar w);u^\perp \ra = -O_u(n)+\Theta_u(n^3)=\Theta_u(n^3) \]
for all elements $w\in\Lc'$ such that $\hat w=\hat \gamma_{u,n}$. In particular, $\gamma_{u,n}\notin\ev(\Lc')$ for $n$ large.

\section*{Acknowledgement}

\begin{adjustwidth}{3mm}{3mm}
	I would like to thank Laura Ciobanu for suggesting to look at the Rational Subset Membership in $H_3(\Z)$, Alex Levine for asking about constrained equations, and Tatiana Nagnibeda for her constant support. This article has also greatly benefited from the anonymous referee's exhaustive suggestions and comments.
\end{adjustwidth}

\bibliographystyle{plain}
\bibliography{references}

\end{document}